\documentclass[]{article}
\usepackage[a4paper,margin=3cm]{geometry}
\usepackage[draft]{changes}
\definechangesauthor[name=Jose, color=blue]{JP}
\definechangesauthor[name=Philip, color=red]{PH}
\definechangesauthor[name=Martin, color=purple]{MS}
\usepackage{authblk}
\usepackage{graphicx}
\usepackage{graphics}
\usepackage{siunitx} 
\usepackage{subcaption}
\usepackage{amsmath}
\usepackage{cleveref}
\usepackage{booktabs}
\usepackage{algorithm}
\usepackage{algpseudocode}
\usepackage{amsfonts}
\usepackage{amsopn}
\usepackage{todonotes}
\usepackage{arydshln}
\usepackage{xcolor}
\usepackage{tikz}

\usepackage[giveninits=true, style=numeric,maxbibnames=6, url=true, doi=true, hyperref=false]{biblatex}
\title{Shape optimization in $\InfBanach$ with geometric constraints:\\ a study in distributed-memory systems}

\author[1]{Philip J. Herbert}	
\affil[1]{Department of Mathematics, University of Sussex, Brighton, BN1 9RF, United Kingdom}
\author[2]{Jose A. Pinzon Escobar}
\author[2]{Martin Siebenborn}
\affil[2]{{Department of Mathematics, University Hamburg, Bundesstr. 55, 20146 Hamburg, Germany}}

\DeclareMathOperator{\argmin}{arg\,min}
\DeclareMathOperator{\Div}{div}

\newcommand{\Gin}{{\Gamma_\mathrm{in}}}
\newcommand{\Gout}{{\Gamma_\mathrm{out}}}
\newcommand{\Gwall}{{\Gamma_\mathrm{wall}}}
\newcommand{\Gobs}{{\Gamma_\mathrm{obs}}}
\newcommand{\R}{\mathbb{R}}
\newcommand{\mult}[1]{{\psi_{#1}}}

\newcommand{\InfBanach}{W^{1,\infty}}
\newcommand{\Hilbert}{H^{1}}
\newcommand{\PBanach}{W^{1,p}}
\newcommand{\id}{{\rm id}}

\newcommand{\test}[1]{{\delta_{#1}}}
\newcommand{\stest}[1]{{\mu_{#1}}}

\newcommand{\state}{y}
\newcommand{\visc}{\nu}
\newcommand{\press}{{\mathfrak{p}}}
\newcommand{\vel}{{{v}}}
\newcommand{\defor}{u}

\newcommand{\gmult}{\mult{g}}

\usetikzlibrary{backgrounds,patterns,calc}

\newcommand{\holdall}{{\rm D}}

\begin{document}
\maketitle
\begin{abstract}
{In this paper we present a shape optimization scheme which utilizes the alternating direction method of multipliers (ADMM) to approximate a direction of steepest descent in $\InfBanach$.
The followed strategy is a combination of the approaches presented in Deckelnick, Herbert, and Hinze, ESAIM: COCV 28 (2022) and M\"uller et al.~SIAM SISC 45 (2023).
This has appeared previously for relatively simple elliptic PDEs with geometric constraints which were handled using an ad-hoc projection.
Here, however, the optimization problem is expanded to include geometric constraints, which are systematically fulfilled.
Moreover, this results in a nonlinear system of equations, which is challenging from a computational perspective. 
Simulations of a fluid dynamics case study are carried out to benchmark the novel method.
Results are given to show that, compared to other methods, the proposed methodology allows for larger deformations without affecting the convergence of the used numerical methods. 
The mesh quality is studied across the surface of the optimized obstacle, and is further compared to previous approaches which used descents in $\PBanach$. 
The parallel scalability is tested on a distributed-memory system to illustrate the potential of the proposed techniques in a more complex, industrial setting.
}

\par

\end{abstract}

\section{Introduction}
In this paper we investigate the efficient optimization of shapes where the optimization must take geometrical constraints into consideration.
We are particularly interested in so-called PDE-constrained shape optimization, where the objective functional, which depends on a domain, also depends on the solution of a PDE within that domain.
Shape optimization is well studied and we refer to \cite{sokolowski1992,zolesio2011,haslinger2003introduction} for an overview.
Here, special interest is placed in the case where only a subset of the boundary should be deformed, the so-called obstacle, and where it should be possible to experience large deformations with minimal degradation of the mesh.
The example considered here is that of an object, which we will refer to as an obstacle, inside a flow tunnel.
Inside this flow tunnel the flow is described by a PDE, such as Stokes or Navier--Stokes, and the surface of the obstacle is optimized for a given functional.
We will consider only the case where the PDE constraint is given by the non-linear stationary Navier--Stokes equations and the functional is the energy dissipation.
In this setting, the obstacle should maintain a fixed volume and barycenter, otherwise the minimizer would be no obstacle.
There is widespread interest in this topic, for instance \cite{schmidt2013three,berggren2016largescale} and the algorithmic approaches considered in \cite{mueller2021,siebenborn2018shape}.
This problem is classically studied  and the solution is described in \cite{pironneau1973optimum,pironneau1974optimum} as being a prolate-spheroid or rugby ball, given that tips have to formed on the surface of the object as part of the optimization process. 
Additionally, this problem was explored more recently in a shape optimization context in \cite{mohammadi2010,blauth2021nnl,schulz2016}.

In shape calculus, one is often interested in using a descent method to obtain a series of shapes which should approximate a minimizer of the objective, $J$.
This requires the computation of the first derivative of the objective with respect to the domain, which we will refer to as the shape derivative, $J'$.
A widely used approach to finding a descent direction is based on Hilbert spaces.
In such an approach, one relates the shape derivative $J'$ to the so-called shape gradient $\nabla J$ by use of Riesz representation in the topology of the Hilbert space.
This is the case, for instance, in \cite{vogel21,pinzon2020parallel}.
A variety of this is the so-called extension equation approach, \cite{onyshkevych2020,mueller2021,haubner2021,vogel21}, where an elliptic PDE is solved to find a deformation field in $\Hilbert$.
For further details, we refer to the overview article \cite{allaire2020} and the bibliography within.

Recent trends have involved the use of Lipschitz, or $\InfBanach$, functions and their approximation by $\PBanach$ functions.
The case of $\InfBanach$, to the best of our knowledge, was first introduced in a practical setting in \cite{deckelnick2021} where it was restricted to the optimization of star-shaped domains.
Generic, non-star shaped domains were considered in \cite{deckelnick23} using deformations in $\InfBanach$.
The approximation with functions in $\PBanach$ is based on \emph{relaxing} the space and choosing $p$ \emph{large enough} so that one is, in some sense, close to being a minimizer in $\InfBanach$.
A method which handled geometric constraints by using an augmented Lagrangian method was introduced in \cite{mueller2021}.
In order to handle \emph{large} values of $p$, an iterative approach which incrementally increased $p$ was used.
The convergence, optimal achievable shape, and mesh quality were compared to a Hilbert space method approach. 
In \cite{mueller2022scalable}, an algorithm for handling geometric constraints of the integral type is combined with the $p$-Laplace relaxation to find descent directions that preserve volume and barycenter.  

For many applications, the initial configuration must go through large deformations to reach an optimal shape.
While, mathematically and in the real world, there is no mesh, for computational methods, it is often utilized.
Large deformations can lead to a loss of mesh quality, i.e.~degeneration of discrete grid elements. 
Degeneration of elements has a negative impact on the convergence of iterative solvers as well as the approximation properties of the discrete solution.
It is therefore important to take into account the mesh quality, if remeshing is to be prevented. 
Preservation mechanisms are often taken into account in computational methods.
In \cite{onyshkevych2020,pinzon2021fluid,haubner2021}, a constraint on the determinant of the deformation gradient was imposed, and enforced using an augmented Lagrangian.
The incorporation of such constraints may reduce the space of attainable shapes.

This article firstly wishes to build upon the work of \cite{mueller2022scalable}, which considers computational scaling to utilize the steepest descent methods in $\InfBanach$ discussed in \cite{deckelnick23}.
Secondly, we investigate, to a limited extent, the quality of the produced mesh.
Finally we demonstrate the computational scalability of the method, which is in no way guaranteed from the previous work which studied $W^{1,p}$ as the scheme is entirely different.

The article is structured as follows:
In Section \ref{sec:model}, we present the essential background for shape optimization in $\InfBanach$ \cite{deckelnick23}, as well as the benchmark physical model which will be the subject of study.
Section \ref{sec:optim} introduces the algorithms used within the novel components in this work.
Simulation results for the $\InfBanach$ methodology appear in Section \ref{sec:results}.
We present in Section \ref{sec:comp} a comparison between the resulting meshes and objective function for $\InfBanach$ and $\PBanach$ approaches. 
Finally, the parallel scalability for the $\InfBanach$ method is measured in Section \ref{sec:hpc}.

\section{Shape optimization in $\InfBanach$ and the model problem}\label{sec:model}
This and the following sections build heavily on the work presented in \cite{mueller2022scalable,deckelnick2021,deckelnick23}.
An in-depth discussion of the theoretical aspects of shape optimizations is not within the scope of this work, therefore interested readers are referred to well-known monographs such as \cite{sokolowski1992,zolesio2011,henrot2018}.
It is however necessary to recount a few aspects in order to build towards the methodology described in \Cref{sec:optim}.

\subsection{Shape Optimization}\label{ss:shapeopt}
The task of (PDE constrained) shape optimization is to, given a collection of admissible domains $\mathcal{S}_{ad}$ and a functional on those domains $J\colon \mathcal{S}_{ad}\to \R$, find $\Omega^* \in \mathcal{S}_{ad}$ which attains the minimum.
Typically $\mathcal{S}_{ad}$ will be a collection of open and bounded domains domains in $\R^{d}$ for $d \in \{2,3\}$.
It may also be the case that one wishes to restrict the minimization to be over domains which satisfy a (geometric) constraint, say $g(\Omega) = 0$.

Finding minimizers in practice is a difficult task.
Iterative methods are often used to find stationary points.
A standard strategy is to, given the current domain $\Omega$, update the domain $\Omega_{new} := \Omega(\defor):= \{ x + \defor(x) : x \in \Omega\}$ for some suitably chosen $\defor$.
For convenience, we denote the perturbation of the identity by $F:= \id + u$, so that $\Omega_{new} = F(\Omega)$.
This means that our iterates yield a natural parameterization over the initial domain.

The task of choosing $\defor$ in practice is non-trivial and many such choices can be made, as discussed in the introduction.
In this work, we will consider the steepest descent in $\InfBanach$.
Where possible, parallels are drawn to the $p$-Laplace optimization scheme proposed in \cite{mueller2021,mueller2022scalable}.
A particular reason for using $\InfBanach$ functions is the following fact: When $\Omega$ is convex and $|D \defor| < 1$ a.e.~in a sub-multiplicative norm, then the map $F$ is bi-Lipchitz.
Where $D \defor$ is the Jacobian matrix of $\defor$. 
Two examples of sub-multiplicative norms are the spectral and Frobenius norms.
In the case that $\Omega$ is not convex, this bi-Lipschitz property need not hold, however one has that for any convex subset $K$ of $\Omega$, the restriction $F|_{K}$ is bi-Lipschitz onto its image.
In \cite{deckelnick23} it is suggested to account for potential non-convexity by the use of a fictitious domain.
For the situation we consider, this is not done.
In the experiments the constructed maps remain bi-Lipschitz.

To formulate a direction of steepest descent, it is useful to have a derivative.
We say that $J$ is shape differentiable at $ \Omega$ if the map $\defor \mapsto J(\Omega(\defor))$ is suitably differentiable at $u=0$ in $\InfBanach$.
We denote the directional derivative as $J'(\Omega)u:=\lim_{t\to 0}\frac{J(\Omega(t\,u))\,-\,J(\Omega)}{t}$.
For the problem we consider, the derivative $J'$ is well known, however one may wish to use the so-called Lagrange multiplier method to calculate it, see \cite{hinzeulbrich2009} for example.
In order to handle the geometric constraints, we will make use of a non-linear space to find the steepest descent so that at each iteration the constraints are fulfilled to a given tolerance rather than approximated, which is the case when using, for example, an augmented Lagrange approach.
This direction of descent is found as: given $\sigma \in (0,1)$ find
\begin{equation}\label{eq:geometricDirectionOfDescent}
\begin{aligned}
	\min_{\defor \in \InfBanach(\Omega;\R^d)} J'(\Omega) \defor
	\\
	\mbox{s.t. } \defor|_{\partial \holdall }=0,
	\\
	\|D \defor \|_{L^\infty(\Omega)} \leq \sigma,
	\\
	g(\Omega(\defor)) = 0,
\end{aligned}
\end{equation}
where $\holdall$ is a sufficiently smooth hold-all domain which will be introduced later.
It is worth mentioning that due the geometric constraints, $J'(\Omega)\defor$ need not equal $-\|J'(\Omega)\|$ as in e.g.~\cite{deckelnick23}, where
\begin{equation}
	\|J'(\Omega)\| := \sup \{ J'(\Omega)\tilde \defor : \tilde \defor \in \InfBanach(\Omega;\R^d),\, \tilde \defor|_{\partial \holdall }=0,\, \|D \tilde \defor \|_{L^\infty(\Omega)} \leq 1,\, g'(\Omega)\tilde \defor  = 0 \}
\end{equation}
is the (dual) norm of $J'(\Omega)$ on the subspace which corresponds to constraining $g$.
One may however expect that $J'(\Omega)\defor / \sigma$ converges to $-\|J'(\Omega)\|$ as $\sigma \to 0$.
Such an expectation comes from the notion of the Hadamard derivative.
The existence of such a $\defor$ in the continuous has not yet been developed, but will be considered in upcoming work.
Let us note that obtaining a function satisfying \eqref{eq:geometricDirectionOfDescent} is not necessarily trivial; the construction we specifically use in the experiments is given by an ADMM.
Full details are given in Section \ref{subsec:ADMM}


Now that we have discussed the shape optimization framework we are interested in, let us discuss the PDE constraint.

\subsection{Physical model}
We are interested in the case that the energy $J$ is given by
\begin{equation}
J(\Omega) = \frac{\nu}{2} \int_\Omega D v : D v \, dx,
\label{eq:EnergyDispObjectiveFunctional}
\end{equation}
where $\nu>0$ is given and, for a given prescribed inflow $\vel_\infty$, the velocity $v$ weakly satisfies the incompressible stationary Navier--Stokes equations
\begin{equation}
	\begin{aligned}
		- \nu \Delta \vel + (\vel\cdot \nabla)\vel + \nabla \press = 0 & \text{ in } \Omega\\
		\Div{\vel} = 0  & \text{ in } \Omega\\
		\vel = 0  & \text{ on } \Gobs\cup\Gwall\\
		\vel = \vel_\infty  & \text{ on } \Gin\\
		\nu D\vel \,\cdot n = \press n & \text{ on } \Gout.
	\end{aligned}
	\label{eq:NavStokes}
\end{equation}
The regions $\Gobs$, $\Gwall$, $\Gin$, and $\Gout$ are to be described shortly.
In this setting, $\nu$ may be referred to as the kinematic viscosity, and $\vel$, the velocity.
We will assume that the solution pair $(\vel, \press)$ is unique, up to an additive constant for the pressure $\press$.

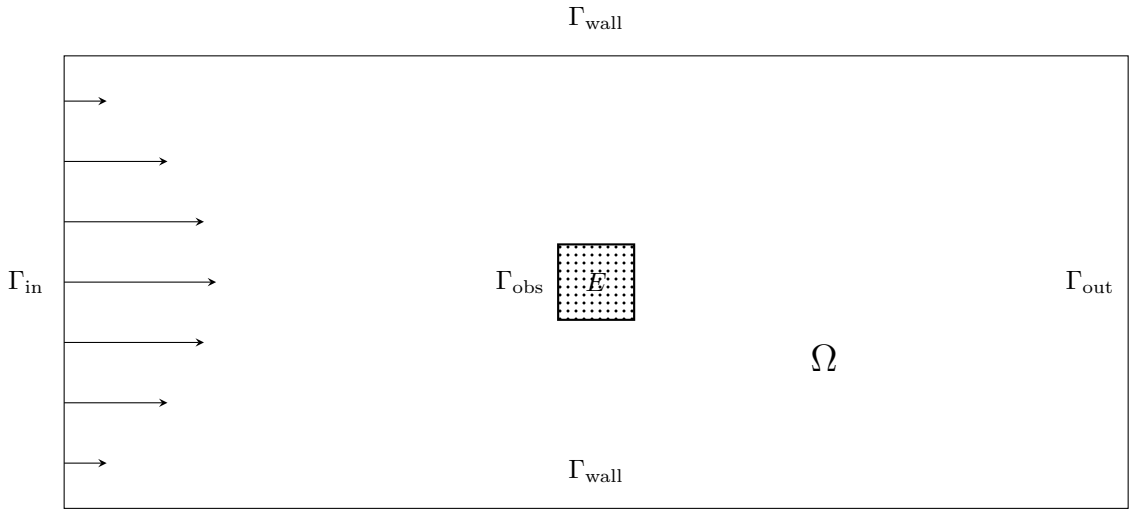
\begin{figure}
	\centering
	\begin{tikzpicture}
		\draw[thick, pattern=dots] (-0.5,-0.5) rectangle (0.5,0.5);
\draw (-7,-3) rectangle (7,3);
\draw (0.0,0.0) node {$E$};
\draw (-7.5,0.0) node {$\Gin$};
\draw (6.5,0.0) node {$\Gout$};
\draw (-1.0,0.0) node {$\Gobs$};
\draw (0.0,3.5) node {$\Gwall$};
\draw (0.0,-2.5) node {$\Gwall$};
\draw (3.0,-1.0) node {\Large $\Omega$};
\foreach \y in {-2.4,-1.6,...,2.4}{
	\draw [thin,-stealth] (-7,\y) -- (-5-0.25*\y*\y,\y);
};
	\end{tikzpicture}
	\caption{A 2d view of the used domain.} \label{fig:domain2d}
\end{figure}

We are interested in the case of an open and bounded domain $\holdall \subset \R^d$, where again $d\in \{2,3\}$, which will take the physical interpretation of a flow tunnel.
Inside this flow tunnel we consider an obstacle, represented by an open and simply connected set $E\subset D$.
The obstacle $E$ should satisfy $ \partial E \cap \partial \holdall = \emptyset$, i.e.~the obstacle should not touch the boundary of the flow tunnel.
We are then interested in optimizing the domain $\Omega := \holdall \setminus \bar{E}$.
For convenience we denote $\Gobs := \partial E$, the boundary of the obstacle.
This has, in words, described the collection $\mathcal{S}_{ad}$.
The full details of the necessary regularity of the domains is beyond the scope of this work.
A sketch of the flow tunnel appears in Figure \ref{fig:domain2d}.

Let us note that the way in which we describe the domain $\Omega$ means that one has that $\partial \Omega = \partial \holdall \cup \Gobs$, where we are not to perturb the flow tunnel $\holdall$, only the obstacle i.e.~$\Gobs$.

Let us now discuss the physical motivation for the constraints which we consider.
It is expected that reducing the volume of $E$ would lead to a lower value for \eqref{eq:EnergyDispObjectiveFunctional}.
For instance, $E$ could vanish, leading to a degenerate minimizer. 
Furthermore, if the obstacle were to move to an area which is less impactful on the overall energy, e.g.~closer to $\Gwall$, where the velocity is vanishing, or downstream within the flow tunnel, one may again expect a reduction in value.
To prevent these cases, geometric constraints are imposed on $\Omega$.
These constraints are on the (normalized) barycenter and volume, and are given by
\begin{align}
	g_i(\Omega(u)) := \int_{\Omega} \left( F_i \det(DF) - \id \right) \; dx &= 0\mbox{ for } i = 1,\ldots,d \label{eq:BarycenterConstraint},\\
	g_{d+1}(\Omega(u)) := \int_{\Omega} \left(\det(DF) - 1 \right) \; d x &= 0 \label{eq:VolumeConstraint},
\end{align}
which represent the difference between the reference configuration, described via the method of mappings, and the deformed domain.
These expressions are used under the assumption that the barycenter is set to the origin, and that the volume remains constant.

\section{Optimization Method}\label{sec:optim}
In this section, the algorithmic approach to finding a minimizing shape is described.
The algorithm we are using is relatively straightforward besides finding descent directions in $\InfBanach$.
The majority of the section is dedicated to the description of how we choose to find this steepest descent.
We explain the methodology used to approximate a solution to the optimization problem formulated in \eqref{eq:geometricDirectionOfDescent}, which relies on the computation of the shape derivative $J'(\Omega)\defor$ to formulate a descent-like method \cite{allaire2020}.
As mentioned before, the proposed methodology uses the ADMM, as in \cite{deckelnick23}, to approximate a deformation field which preserves the mesh topology between shape iterates. 
Nevertheless, it is worth noting that the following optimization scheme is in no way limited to the use of the ADMM, but could be adapted to other approaches to find $u$ in $\InfBanach$. 
The task of finding a potentially better-suited solver is not only interesting, but an ongoing topic of research. 
In the following, the implemented routines are described, starting with the descent-like method and followed by the ADMM in the context of an iterative shape optimization scheme. 

Henceforth, we will deal only with discrete quantities, that is to say, our shape is now represented by the moving mesh with triangulation $\mathcal{T}_h$.
On that mesh we approximate the solution of stationary Navier--Stokes equation by means of a finite element discretization with lowest order Taylor-Hood elements, and the energy is computed using this finite element solution.

\subsection{Shape optimization algorithm}\label{ss:Temp}

The main routine is outlined in \Cref{alg:descent_alg}, where a descent-like method is used to obtain an optimized geometry from a given arbitrarily shaped initial guess, $\Omega^0$.
Within the outer loop of \Cref{alg:descent_alg}, a series of shape iterates is generated by computing a deformation vector $\defor$ and applying it over the nodes of the current discrete domain $\Omega^k$ to propose a domain $\hat{\Omega}$, as described in Section \ref{sec:model}.
The resulting geometry is carried on to the next step if $J(\hat{\Omega})<J({\Omega^k})$, which requires solving the state equation and computing the objective function. 
On the contrary, if the aforementioned condition is not fulfilled then $\sigma$ is reduced as indicated in line 13. 

\begin{algorithm}
	\caption{Shape Optimization-Descent method}
	\label{alg:descent_alg}
	\begin{algorithmic}[1]
		\Require{Initial shape $\Omega^0$,$\epsilon_1$,$\epsilon_2$, $\sigma$,$N$}
		\State $k = 0$
		\State $\defor^k\gets 0$
		\State{$y^k\gets$ Solve state \eqref{eq:NavStokes}}
		\State Compute $J(\Omega)^k$
		\For{k=0,1,\dots,$M$}
		\State{Compute the shape derivative $J'(\Omega^k)$ using the adjoint method}
		\State{Modification of $J'(\Omega^k)$ as in Section \ref{ss:ShapeDer}}
		\State{$\defor^k\gets\textsc{ADMM}(J'(\Omega),\sigma,\epsilon_{\mathrm{2}},\epsilon_{\mathrm{3}},N)$} as in Section \ref{subsec:ADMM}
		\State{Temporary geometry $\hat{\Omega}=\Omega^k(u^k)$}
		\State{$y\gets$ Solve state equation $e$ on $\hat{\Omega}$}		\State{Compute objective function $J(\hat{\Omega})$}
		\If{$J(\Omega^k) \leq  J(\hat{\Omega})$}
		\State $\sigma\gets \sigma/2$
		\Else
		\State $y^{k+1}\gets y$
		\State $\Omega^{k+1}\gets\hat{\Omega}$
		\State $k= k + 1$
		\EndIf
		\EndFor
	\end{algorithmic}
\end{algorithm}

The core of \Cref{alg:descent_alg} lies on finding the descent direction $\defor$, shown in line 8.
This requires the solutions to the state equation \eqref{eq:NavStokes} and its adjoint, for the computation of the shape derivative $J'(\Omega)$.
The solution $\state$ is used to compute the objective function on the temporary geometry, which determines if a descent direction has been found.

The deformation field constitutes a descent direction such that, for the solution of the equation system \eqref{eq:optimality_system_defor}, $J'(\Omega)u<0$ on every step. 
If one were to have unlimited computational time, one might choose to use the descent as the convergence criterion.
Where we recall that $J'(\Omega) \defor$, when rescaled, approximates the directional derivative.
In such a situation, for some $\epsilon_1>0$, one might choose to run the outer loop until 
\begin{equation}
	{|-J'(\Omega^k)\defor^k|} \,<\, {\sigma} \epsilon_{\mathrm{1}},
	\label{eq:ConvCriterion}
\end{equation}
which would, for $\sigma$ small, roughly correspond to $\|J'(\Omega^k)\| < \epsilon_1$.

\subsection{Modification of the shape derivative in the interior of the mesh}\label{ss:ShapeDer}
It is well known by the Hadamard structure theorem that for sufficiently regular domains, the shape derivative $J'(\Omega)$ should be supported only at the boundary, e.g.~\cite[Section 3.4, Theorem 3.6]{zolesio2011}, moreover that it should be a measure. 
However, it is not the case that the shape derivative evaluated on the solution of the discrete state equation is supported at the boundary only.
This may be attributed to the finite element approximation errors - see \cite{hiptmair2014shape} who discuss the error due to different forms of the shape derivative.
It is suggested in \cite{siebenborn2017algorithmic} that one can remove the values which have contributions only in the interior.
Such an approach has been noted to produce qualitatively better results \cite{schuzsieb2016}.
Computationally, this corresponds to only assembling (integrating) the shape derivative over cells which have a vertex which intersects the boundary $\Gobs$.
This is performed in line 7 of \Cref{alg:descent_alg}.

Mathematically speaking, since we are expecting the shape derivative to be a measure, we may represent this by the discrete functional which is a sum of Dirac deltas over the vertices of the mesh.
The aforementioned expert knowledge then says that the relevant nodes are those which are at the boundary of the obstacle.
This recovers the method of \cite{siebenborn2017algorithmic}.

\subsection{The ADMM for the solution of the steepest descent}\label{subsec:ADMM}
The minimization problem for the direction of steepest descent in \eqref{eq:geometricDirectionOfDescent} is highly non-trivial, both in terms of the constraint on the Jacobian, but also the geometric constraints.
To handle the Jacobian condition, inspired by \cite{bartels2021}, \cite{deckelnick23} utilized the ADMM to find directions of steepest descent for PDE-constrained shape optimization problems.
This approach is here extended to include geometric constraints, building on the scheme presented in \cite{mueller2022scalable}.
In there, a methodology to incorporate the geometric constraints was proposed with the intention of avoiding an augmented Lagrangian.

The function which we will consider for the problem in \eqref{eq:geometricDirectionOfDescent} is given as
\begin{equation}
	L(u,\mult{g},q,\lambda) := J'(\Omega)\, u + \sum_{i = 1}^{d+1} \mult{g,i} g_i(\Omega(\defor))
	+ \frac{\tau}{2} \| Du - q \|_{L^2}^2 + \tau (Du - q, \lambda)_{L^2},
	\label{eq:admmLagrangian}
\end{equation}
where the last two terms appeared in \cite{deckelnick23} as part of the ADMM augmented Lagrangian, cf.~\cite{bartels2020,bartels2021}. 
Here, $\mult{g,i}$ are components of the finite dimensional vector $\mult{g} \in \R^{d+1}$ of Lagrange multipliers associated to the geometric constraints \eqref{eq:VolumeConstraint} and \eqref{eq:BarycenterConstraint}. 
In addition to the expected deformation and Lagrange multipliers, $\defor$ and $\mult{g}$, we also have the appearance of $q$ and $\lambda$ in $L$.
Here $q$ takes the role of a slack variable, and $\lambda$ of a Lagrange multiplier which pushes $q$ towards $D\defor$. 

The ADMM strategy we consider is provided in Algorithm \ref{alg:ADMM}.
It operates by alternatingly minimizing $L$ over $q \in Q_h$ and $\defor \in V_h$ such that the geometric constraints are fulfilled, then a simple update is used for the multiplier $\lambda$. 
The iterative process is repeated until the convergence criteria, line 11 are fulfilled. 
The condition $\lVert\Delta_{\lambda}\rVert^2_{L^2}~+~\lVert\Delta u\rVert^2_{L^2}<\epsilon_2$ is similar to that used in \cite{bartels2021,deckelnick23}, and represents the residual of the updates. 
To ensure that $\lVert Du\rVert_{L^\infty}$ is as close as possible to $\sigma$, we introduce an additional condition on the convergence.
The loop stops when $\lVert D\defor\rVert_{L^\infty}$ is also close to the given $\sigma$, up to a small value $\epsilon_3>0$, which may depend on $\sigma$ itself.

Let us discuss the ad-hoc modification which has been made in lines 12-14 of \Cref{alg:ADMM}.
When $\|J'(\Omega)\|$ is not equal to zero, we expect the deformation field $\defor$ to satisfy $\|D \defor\|_{L^\infty}$ is close to $\sigma$ - in the setting without the geometric constraints, these quantities are equal.
In order to avoid degenerate solutions, the vector which stores the shape derivative $J'(\Omega)$ has its components doubled. 
This could be compared to the variable step-size ADMM presented in \cite{bartels2020}, where the value $\tau$ is determined as part of the algorithm.
Here however, this has not been undertaken to speed up the convergence, but to ensure an appropriate deformation field is found.

\begin{algorithm}
	\caption{Descent Direction in $\InfBanach$}
	\label{alg:ADMM}
	\begin{algorithmic}[1]
		\Require{$J'(\Omega)$,$\sigma$,$\epsilon_2,\epsilon_3,N$}
		\State $\tau\gets 1$
		\State $\defor_0\gets0$
		\For{$i=0,1,...,N$}
		\State Find $q\gets \argmin\{ L(\defor,q,\lambda):q \in Q_h ,|q| \leq \sigma\}$ as in Section \ref{sec:opt_for_q}
		\State Find $\defor\gets \argmin \{L(\defor,q,\lambda): \defor \in V_h,\,  g(\Omega(\defor)) = 0 \}$ as in Section \ref{sec:opt_for_u}
		\State $\Delta_\lambda\gets \tau(D\defor - q)$
		\State $\lambda\gets \lambda + \Delta_\lambda$
		\State $\Delta_\sigma\gets \sigma \,-\, \mathrm{max}(| D\defor|_{L^\infty})$
		\State $\Delta\defor\gets \defor -\defor_0$
		\State $\defor_0\gets \defor$
		\If{($\lVert\Delta_{\lambda}\rVert^2_{L^2}~+~\lVert\Delta u\rVert^2_{L^2}<\epsilon_2$)\texttt{ and }$\Delta_\sigma>\, -\epsilon_3$}
		\If{$\Delta_\sigma>\, \epsilon_3$}
		\State{$J'(\Omega)\gets 2\,J'(\Omega)$}
		\Else
		\State break
		\EndIf
		\EndIf
		\EndFor
	\end{algorithmic}
\end{algorithm}
Within every iteration of \Cref{alg:ADMM}, Newton's method is called to solve the nonlinear optimality system, as shown in line 5.
This is described in Section \ref{sec:opt_for_u}.
The strategy is based on that which is described in \cite[Section 3]{mueller2022scalable}, where the Schur complement operator is computed by using $d+1$ numerical solves.
The computational cost is reduced by using preconditioners with grid-independent convergence, such as the geometric multigrid.

Let us expand further on the minimization with respect to $q$ and the saddle point for $(\defor,\mult{g})$ corresponding to minimization of $\defor$ under the geometric constraints which are utilized for the ADMM algorithm.

\subsubsection{Optimality for $q$}\label{sec:opt_for_q}
For the system involving $q$, it is useful to denote the discrete space
\begin{equation}
	Q_h := \{ q \in L^2(\Omega;\R^{d \times d}) : q|_T \in \mathbb{P}^0(T; \R^{d \times d})\},
\end{equation}
where we note that this is not the only possible discrete space one could choose for $q$, however it retains a level of simplicity.
The optimality system for $q$ is given by
\begin{equation}
	\begin{split}
	\frac{\partial}{\partial q} L (u,\mult{g},q,\lambda)(\test{q} - q) \geq&\, 0 \quad \forall \; \test{q} \in Q_h : |\test{q}|\leq \sigma
	\\
	|q| \leq&\, \sigma
	\end{split}
	\label{eq:optimality_system_q}
\end{equation}
which is a convex minimization problem.
For convenience, we note that
\begin{equation}
	\begin{aligned}
		\frac{\partial}{\partial q} L (\defor, \mult{g},q,\lambda) \test{q}
		=
		\tau \int_{\Omega} (q - D\defor - \lambda):\test{q} \, dx.
	\end{aligned}
	\label{eq:dLdQ}
\end{equation}
This may be decomposed into two steps, solving
\begin{equation}
	\frac{\partial}{\partial q} L (u,\mult{g},q,\lambda)(\test{q} - \tilde q) = 0 \quad \forall \; \test{q} \in Q_h
	\label{eq:optimality_sub-system_q}
\end{equation}
which, observing \eqref{eq:dLdQ}, is seen to correspond to a mass solve, and a pointwise projection
\begin{equation}
	\tilde q \mapsto \frac{\tilde{q}}{\max\left( 1,\frac{|\tilde q|}{\sigma}\right)} =: q.
\end{equation}

\subsubsection{Optimality for $(\defor,\mult{g})$}\label{sec:opt_for_u}
For the optimality system involving $(\defor,\mult{g})$, it is convenient to introduce
\begin{equation}
	V_h := \{ \defor \in C^0(\bar{\Omega};\R^d)  :   u|_{T} \in \mathbb{P}^1(T;\R^d),\, T \in \mathcal{T}_h,\ u|_{\partial D} = 0\}
\end{equation}
which is the natural space for deforming a mesh which has flat triangles.
The optimality system for the saddle point for $(\defor, \mult{g})$ is
\begin{equation}
	\begin{aligned}
		\frac{\partial}{\partial u} L (u,\mult{g},q,\lambda) \test{u} &= 0  &&\quad \forall \; \test{u} \in V_h\\
		\frac{\partial}{\partial \mult{g}} L (u,\mult{g},q,\lambda) \test{\mult{g}} &= 0 &&\quad \forall \; \test{\mult{g}} \in \R^{d+1},
	\end{aligned}
	\label{eq:optimality_system_defor}
\end{equation}
which corresponds to a (discrete) Poisson problem with non-linear constraints.
Following the strategy proposed in \cite{mueller2022scalable}, the system is solved via its linearization in the context of Newton's method.
A particular advantage in our setting is that the second derivative of $L$ with respect to $\defor$ is more convenient to handle:
\begin{equation}
	\frac{\partial^2}{\partial u^2}L (u,\mult{g},q,\lambda) (\test{u}, \stest{u}) = \tau\int_{\Omega}D\test{\defor}:D\stest{\defor}
	+ \sum_{i=1}^{d+1} \langle g_{i} ''(F(\Omega)) \test{u}, \stest{u}\rangle.
	\label{eq:Luu}
\end{equation}
The advantage becomes evident by contrasting \eqref{eq:Luu} to the first integral term of the same derivative presented in \cite{mueller2022scalable}, the second variation of the $p$-Laplace energy, which is a degenerate elliptic operator.
In the $p$-Laplace case, $L_{uu}$ was modified to prevent divide-by-zero operations.
Here, this is not necessary.

To solve the saddle point system which appears in \eqref{eq:optimality_system_defor}, the increments due to the Newton method may be expressed as
\begin{equation*}
	\begin{pmatrix}
		A & B \\
		0 & -B^T A^{-1}B
	\end{pmatrix}
	\begin{pmatrix}
		\test{u} \\
		\test{\lambda} 
	\end{pmatrix} = 
	\begin{pmatrix}
		r_u \\
		r_\lambda - B^T A^{-1} r_u
	\end{pmatrix},
	\label{eq:SchurComplementReducedSystem}
\end{equation*}
where, for readability purposes, we write
\begin{equation}
	\begin{aligned}
		A\test{u} &:= \frac{\partial^2}{\partial u^2} L(u^k, \gmult^k, q^k, \lambda^k)(\cdot, \test{u}) \\
		B\test{\gmult} &:= \frac{\partial}{\partial u\,d\gmult} L(u^k, \gmult^k, q^k, \lambda^k)(\cdot, \test{\gmult})\\
		B^T\test{u} &:=\frac{\partial}{\partial \gmult\,\partial u} L(u^k, \gmult^k, q^k, \lambda^k)(\test{u}, \cdot )\\
		r_u &:= - \frac{\partial}{\partial u} L(u^k, \gmult^k, q^k, \lambda^k) \\
		r_\gmult &:= -\frac{\partial}{\partial \gmult} L(u^k, \gmult^k, q^k, \lambda^k).
	\end{aligned}
	\label{eq:discreteDifferentialOperators}
\end{equation}
Let us comment that the Schur complement operator, $S:=-B^T A^{-1}B$, is related to the computation of the $\mult{g}$ Lagrange multipliers.

\section{Results}\label{sec:results}
We present simulation results for 2d and 3d fluid dynamics case studies, where the domain is as described in \Cref{fig:domain2d}.
The initial obstacle $E$ is set as a square, or a box, in 2d and 3d, respectively.
This showcases the successful creation and removal of geometric singularities through the optimization process.
Results are shown to highlight the large deformations present, particularly at the initial steps, in 2d.
This can be compared to the results presented in similar studies, e.g.~\cite{onyshkevych2020,haubner2021,mueller2022scalable}, to mention some.
Mainly, emphasis is placed on the preservation of mesh quality even under the aforementioned large deformations within a single optimization step. 

These results are generated using UG4~\cite{vogel2014ug4}, a simulation framework tailored for distributed-memory systems~\cite{vogel15,vogel16}.
Most importantly we make use of its geometric multigrid preconditioner \cite{hackbusch94} for the iterative methods used within this work.
Our models were implemented using UG4's plugin functionality, and the performance was measured with its builtin profiler. 
The parallel communication is MPI-based~\cite{mpi40}.
The computational grids were created with GMSH~\cite{gmsh}, using triangular (2d) and tetrahedral (3d) elements.
ParMetis~\cite{parmetis} was used for grid partitioning. 

We recall that the state equation is described by the weak form of the incompressible Navier-Stokes equations and is discretized via the mixed Taylor-Hood finite element, with the lowest order elements $P_2-P_1$.
The adjoint equation is similarly discretized with this.
The viscosity is set to $\visc=0.02$ for all simulations. 
The dimensions of the flow tunnel, the holdall domain, for 2d and 3d are
\begin{equation*}
	\holdall_{2d}=(-7,7)\times(-3,3) \text{ and } \holdall_{3d}=(-10,10)\times(-3,3)\times(-3,3),
\end{equation*}
respectively.
The obstacles are given by
\begin{equation*}
	E_{2d} = (-0.5,0.5)\times (-0.5,0.5) \text{ and } E_{3d} = (-0.5,0.5) \times (-0.5,0.5) \times (-0.5,0.5)
\end{equation*}
and are not triangulated, since we are interested primarily on the optimization of its outer surface $\Gobs$.
At $\Gamma_{in}$ the inflow profile is described by
\begin{equation*}
	\vel_\infty(x) = e_1 \prod_{i=2}^d \cos \left( \frac{\pi}{3} x_i \right) \in \R^d
\end{equation*}
where $e_1:= (1,0,\ldots) \in \R^d$.
Note that $\vel_\infty \cdot e_1 \in [0,1]$ on $\Gamma_{in}$.

We recall that \eqref{eq:optimality_system_defor} is approximated with vector-valued $P_1$ finite elements. 
Following the methodology described in \cite{mueller2022scalable}, the Lagrange multipliers of the geometric constraints, $\gmult$, are not associated to a finite element discretization.
As described in \Cref{sec:optim}, their solution is obtained via a direct solver of an $m\times m$ system of equations, i.e.~the Schur complement system $S$.

The code used for the simulation can be found in \cite{gmgshapeopt2021}.

\subsection{Simulations in 2d}\label{ss:2dresults}
In these a studies, a computational mesh with \num{282624} triangular elements was used.
The surface of the obstacle consists of \num{512} edges. 
We utilize the spectral norm for \eqref{eq:geometricDirectionOfDescent}, in this setting.
In order to reduce the dissipated energy, the deformations must remove the preexisting geometric singularities, which in this case are the corners of $\Gobs$.
At the same time, the necessary tips must be created parallel to the direction of the flow.

\begin{figure}[!htbp]
	\centering
	\includegraphics[width=0.9999\textwidth]{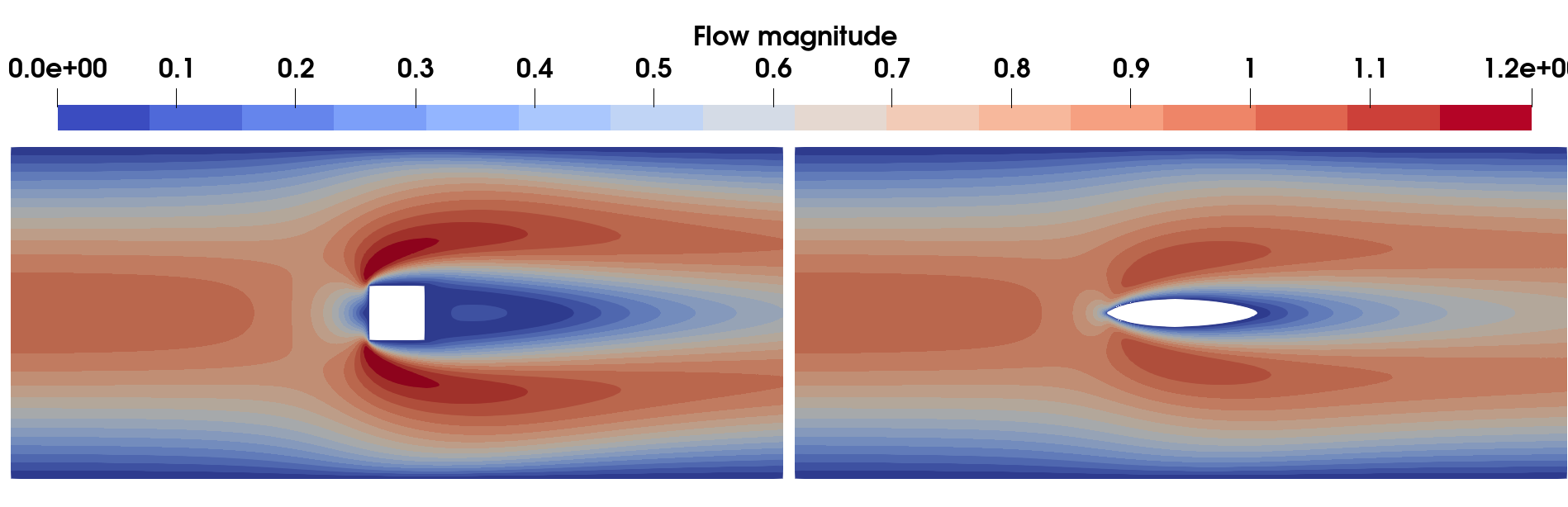}\\
	\caption{2d simulation results for a \num{282624} element grid. The streamlines are shown over the initial (top) and optimized (bottom) $\Gobs$ configurations.}
	\label{fig:2dFinalGrids}
\end{figure}

\begin{figure}[!htbp]
	\centering
	\includegraphics[width=0.7\textwidth]{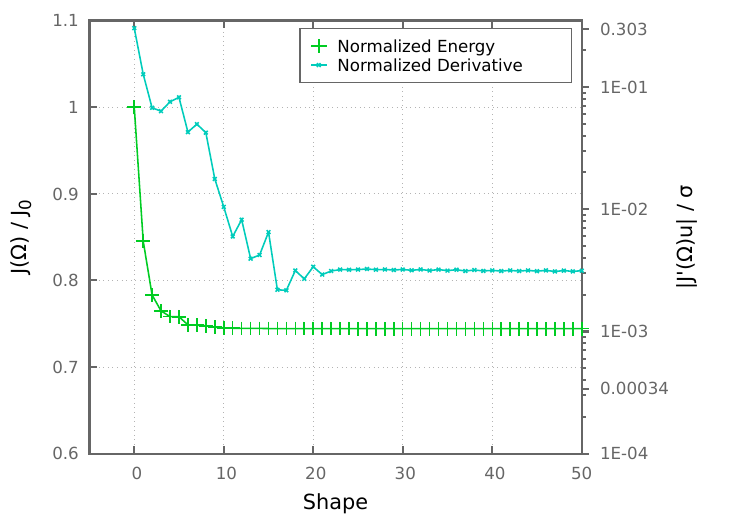}
	\caption{The objective function $J(\Omega)$ divided by the initial value  - $J_0$ - over 50 steps is shown, together with the convergence criterion \eqref{eq:ConvCriterion}.}
	\label{fig:PlotObjFunc}
\end{figure}

The optimization process of \Cref{alg:descent_alg} is illustrated by the results given in \Cref{fig:2dFinalGrids}, where the removal and creation of grid singularities is evident. 
This optimized object surface corresponds to 50 steps of the outer loop of the shape optimization method of \Cref{sec:optim}.
The resulting grid for these simulations appears later in \Cref{sec:comp} and is shown in the upper row of \Cref{fig:FinalGridsComp}.
It shows the regions of the surface $\Gobs$ where singularities were created and removed.
In the figure, the front tip is featured, together with the area corresponding to the upper left corner of the obstacle in the reference configuration.
As previously mentioned, tips must be generated parallel to the direction of the flow and the corners must be smoothed out.
These show that, although the obstacle's surface experiences very large deformations already on the initial deformation steps, the elements retain their shape and, since the $\InfBanach$ deformation is assured to yield a bi-Lipschitz map, no overlapping between elements occurs.

The plot in \Cref{fig:PlotObjFunc} shows the energy dissipation divided by the value calculated on the initial configuration, $J_0:=J(\Omega_0)$.
Together with the objective function, a scaling of the absolute value of the directional derivative is plotted.
Let us recall the discussion after \eqref{eq:ConvCriterion} which comments that this rescaled directional derivative can potentially be related to the norm of the derivative.
The observed monotonic decrease of $J(\Omega)$ is related to the condition in \Cref{alg:descent_alg}, where its checked whether the new geometry $\hat{\Omega}$ obtained from the calculated descent direction reduces $J(\Omega)$.
During the first steps, a very large reduction of the objective function and the derivative can be observed. 
The energy dissipation is reduced in the range of 25\% from the initial value, while the directional derivative is reduced several orders of magnitude, which is to be expected.

\begin{figure}[!htbp]
	\begin{center}
		\begin{tabular}{lccc}
			Step & $\Gamma_{obs}$& Corner & Tip\\
			0 &\includegraphics[width=0.35\textwidth]{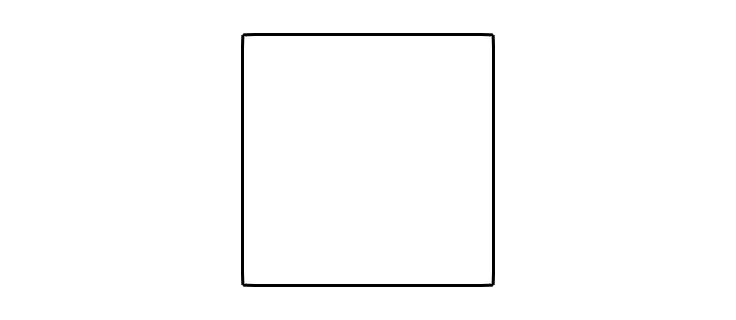}& \includegraphics[width=0.22\textwidth]{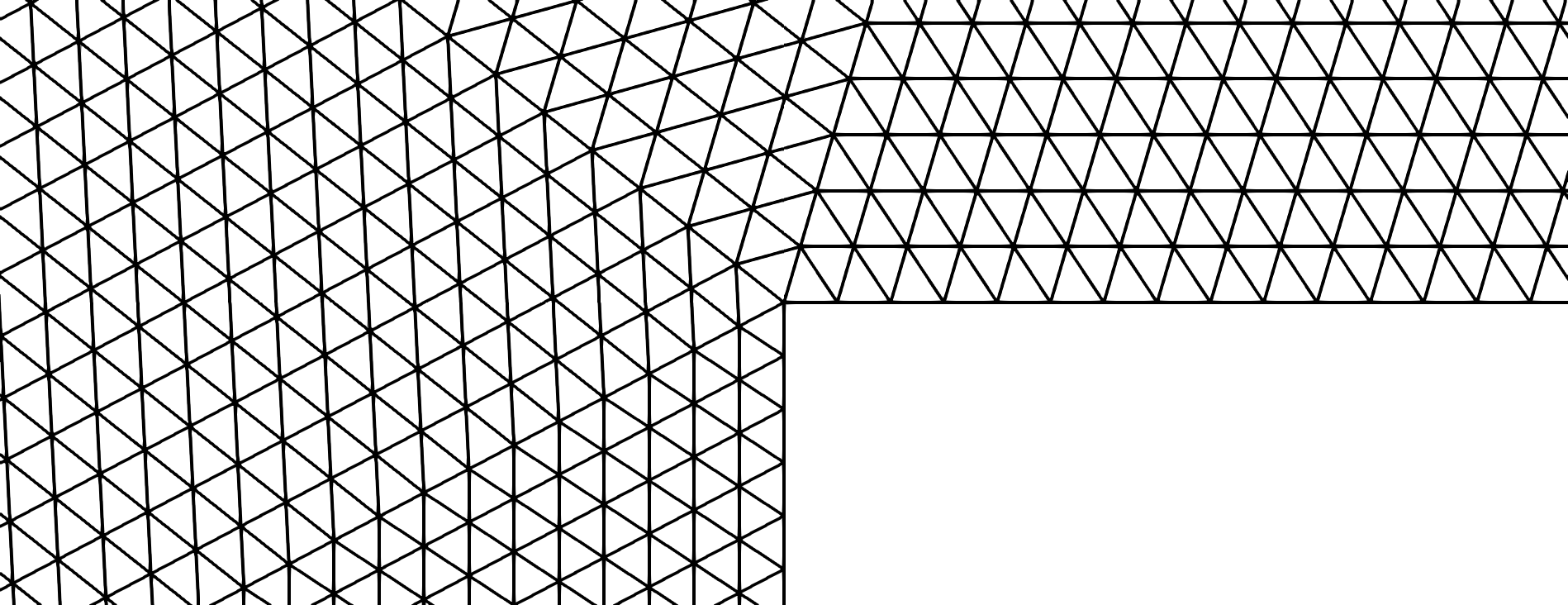}& \includegraphics[width=0.22\textwidth]{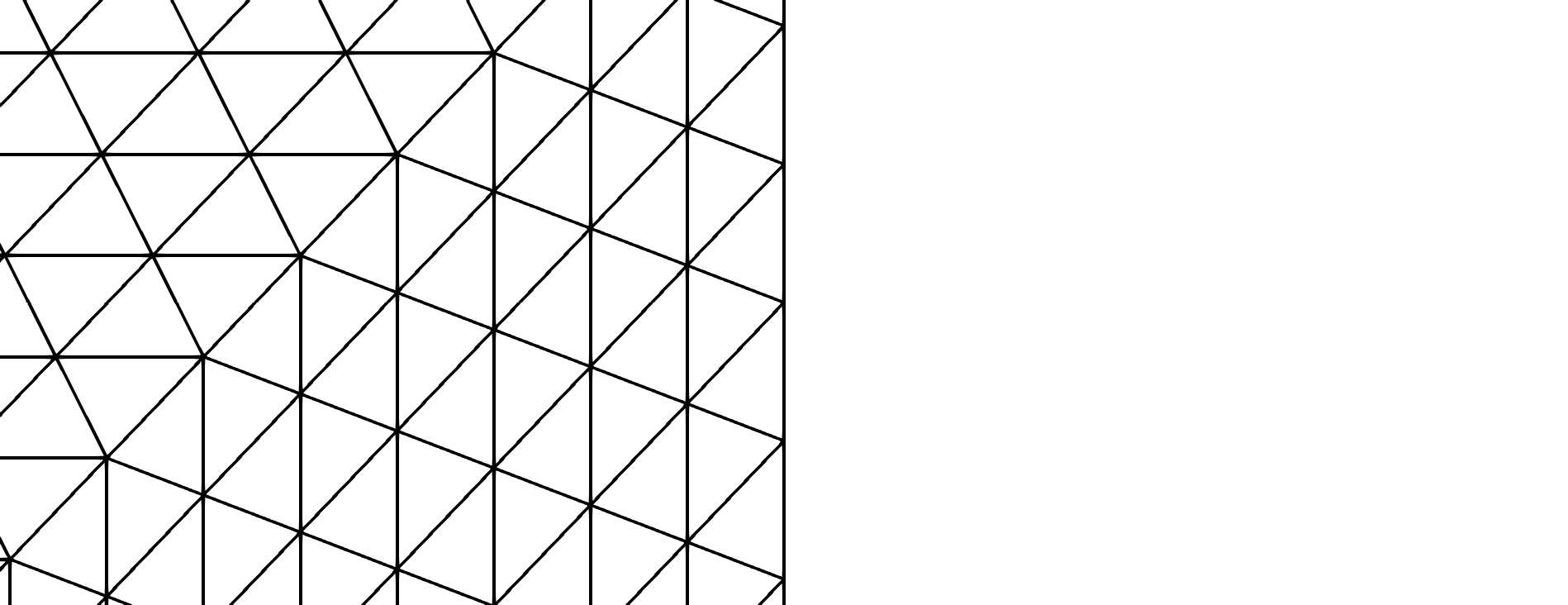}\\
			5 &\includegraphics[width=0.35\textwidth]{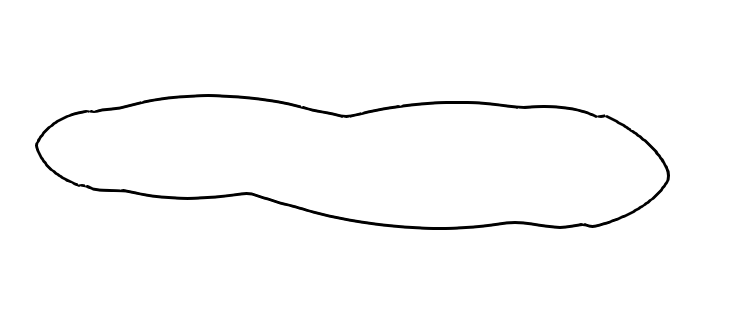}& \includegraphics[width=0.22\textwidth]{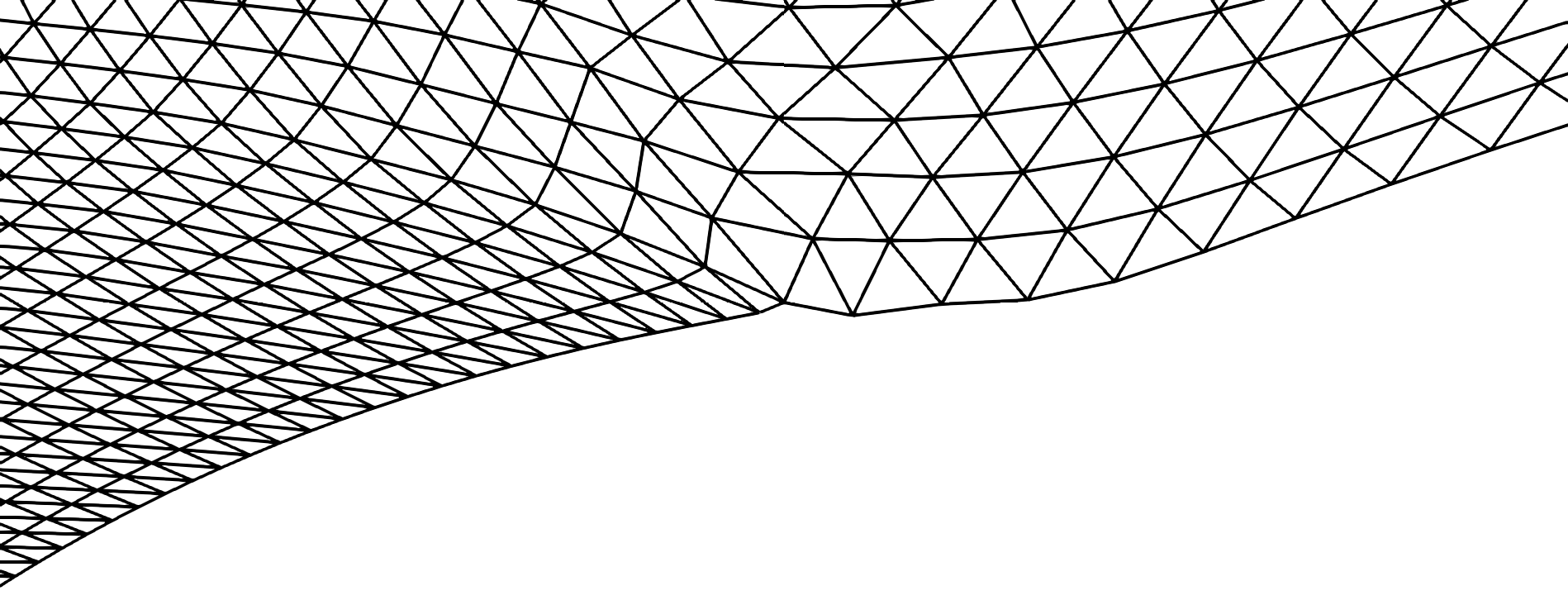}& \includegraphics[width=0.22\textwidth]{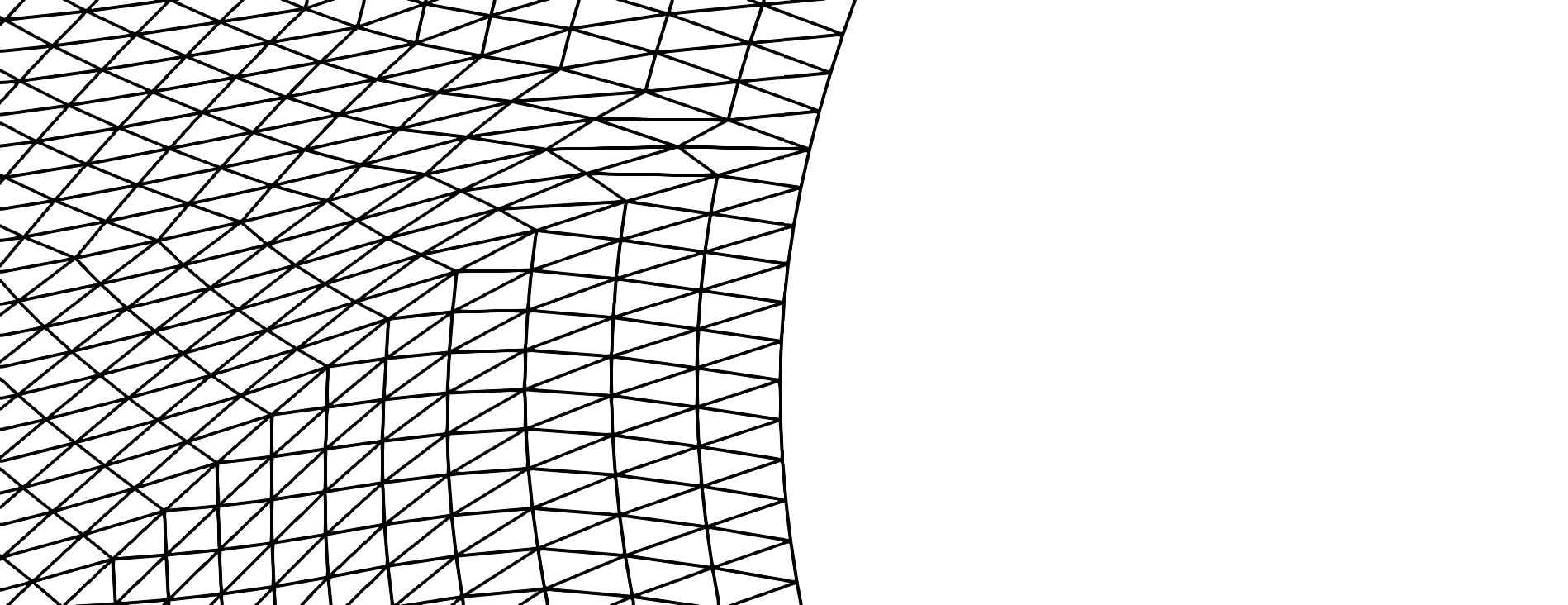}\\
			10 &\includegraphics[width=0.35\textwidth]{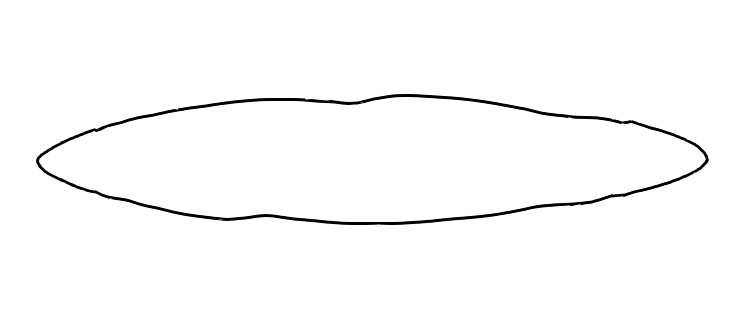}& \includegraphics[width=0.22\textwidth]{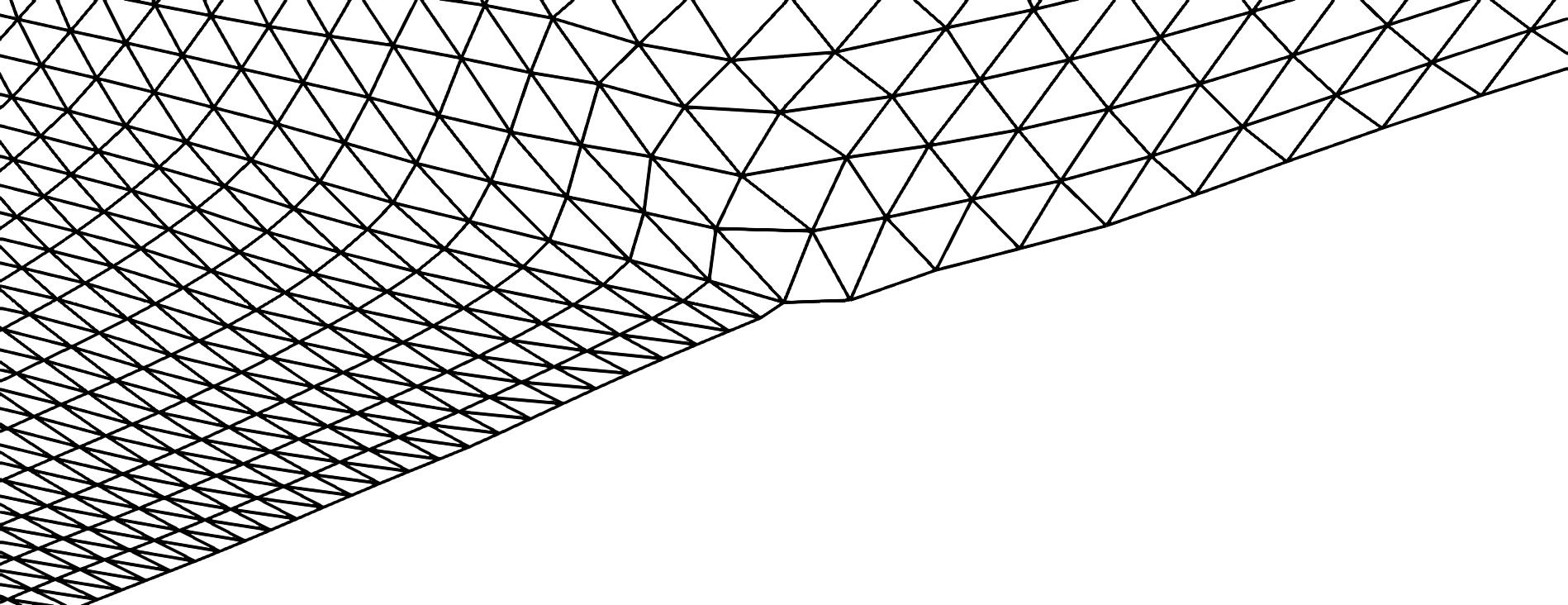}& \includegraphics[width=0.22\textwidth]{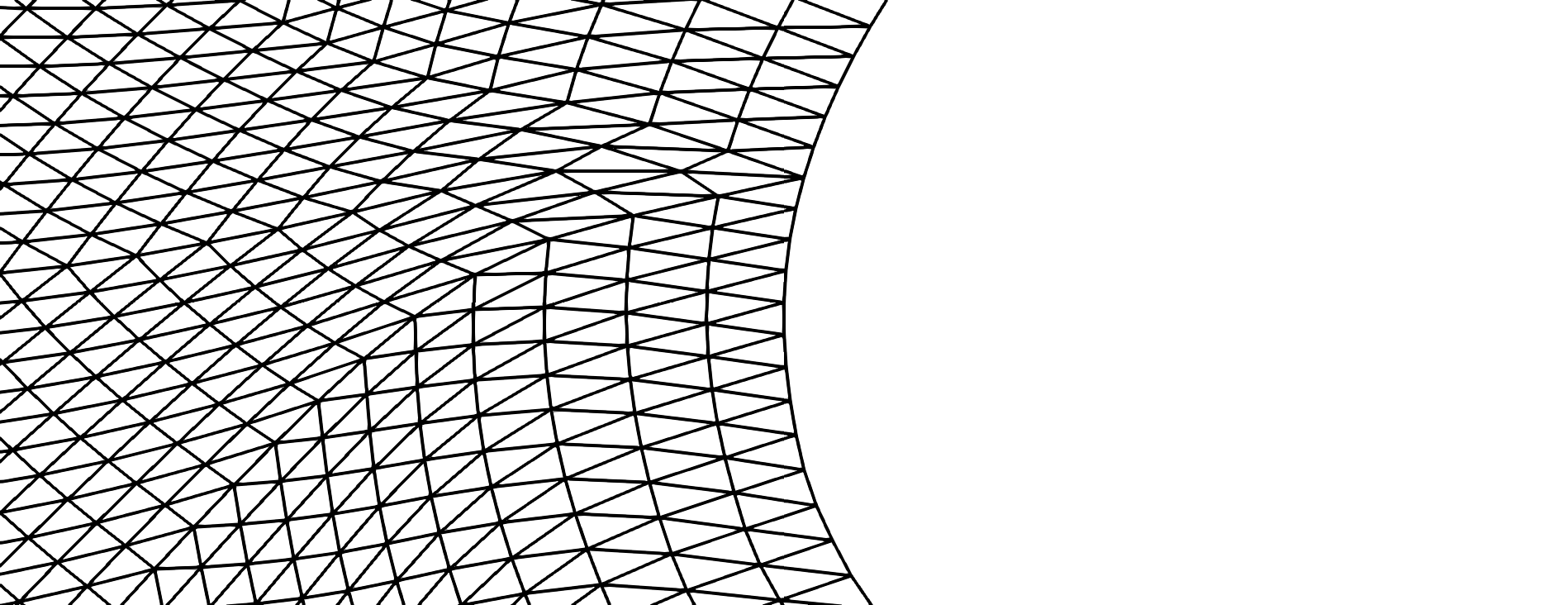}\\
			15 &\includegraphics[width=0.35\textwidth]{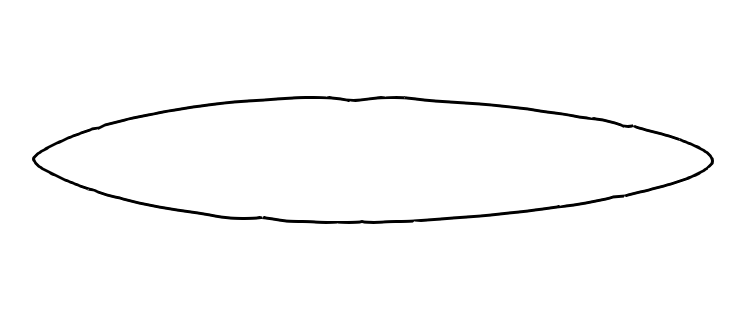}& \includegraphics[width=0.22\textwidth]{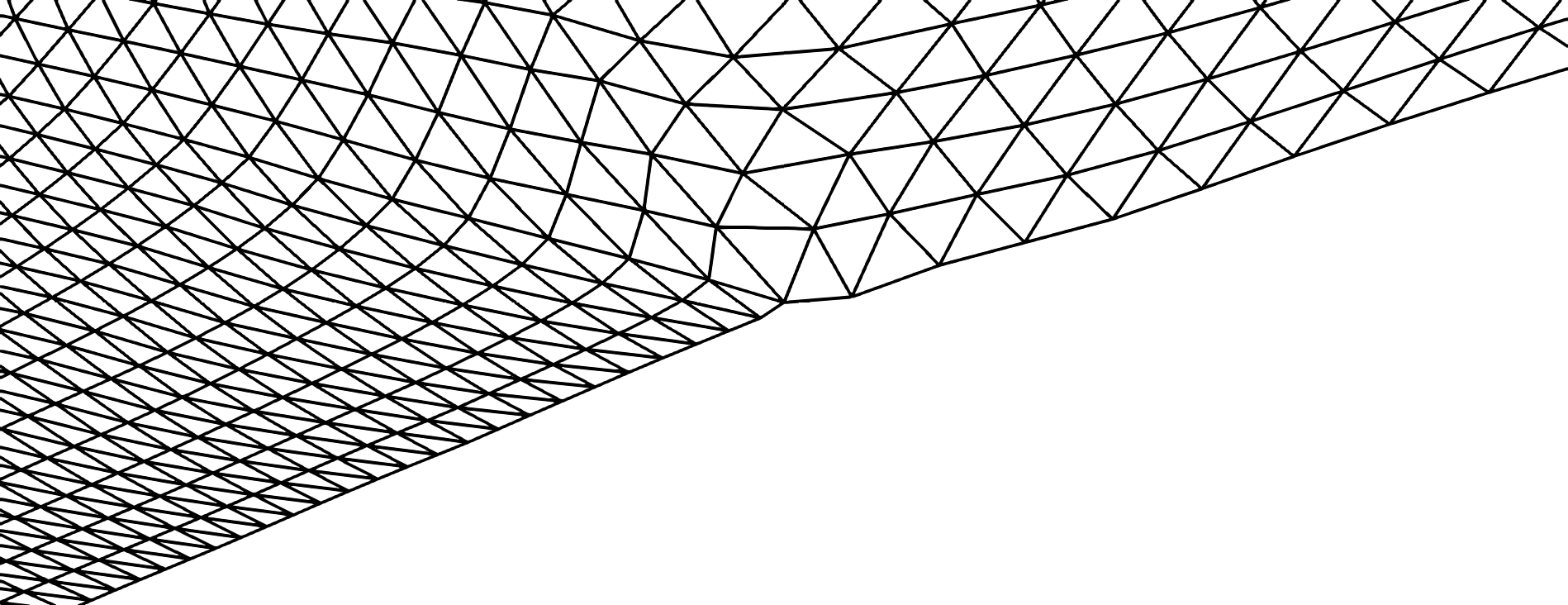}& \includegraphics[width=0.22\textwidth]{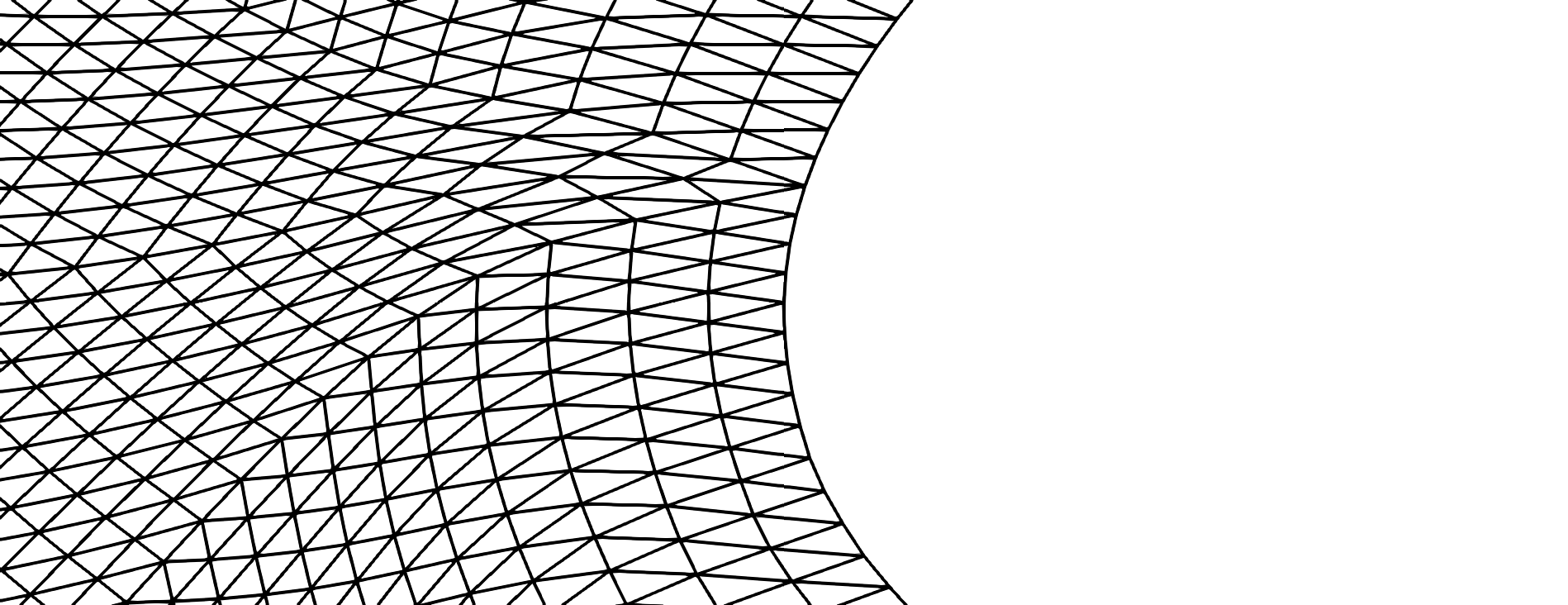}\\
			20 &\includegraphics[width=0.35\textwidth]{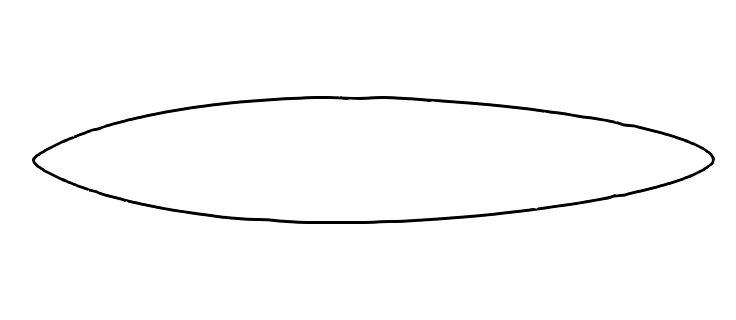}& \includegraphics[width=0.22\textwidth]{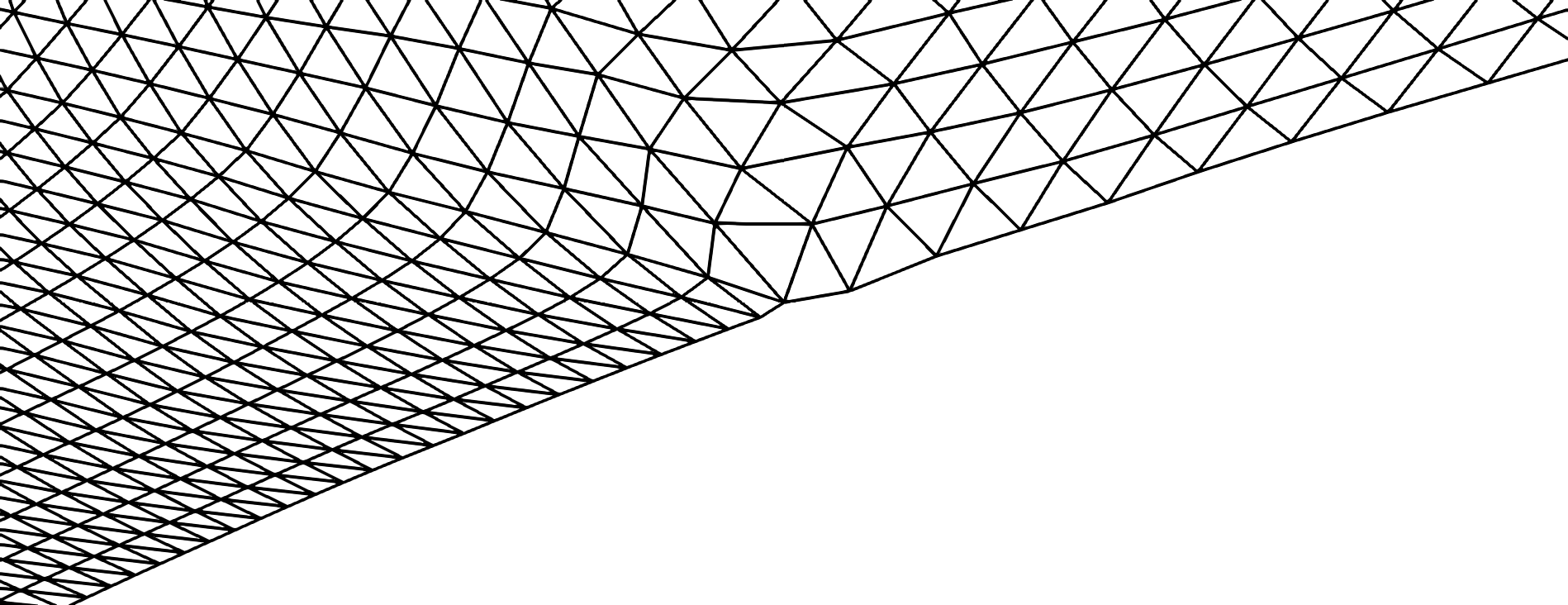}& \includegraphics[width=0.22\textwidth]{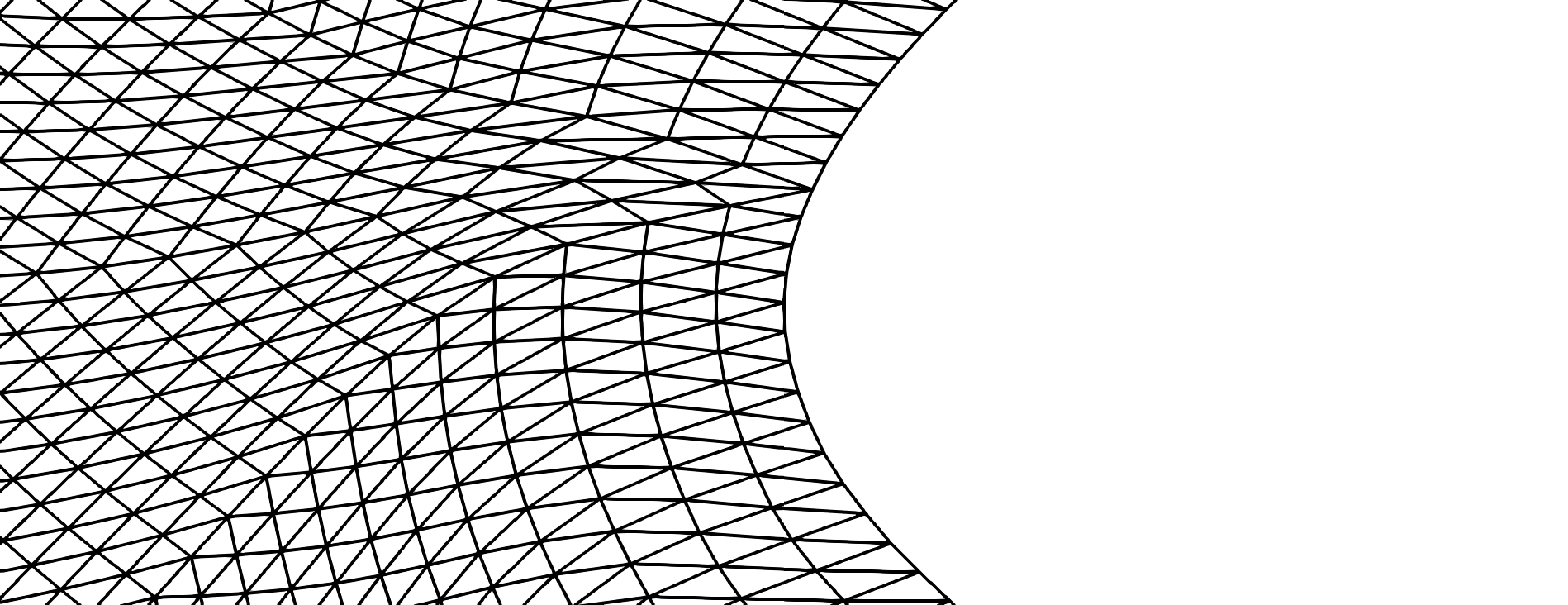}\\
			40 &\includegraphics[width=0.35\textwidth]{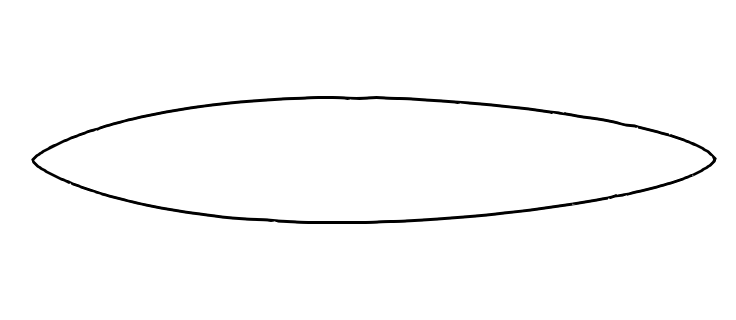}& \includegraphics[width=0.22\textwidth]{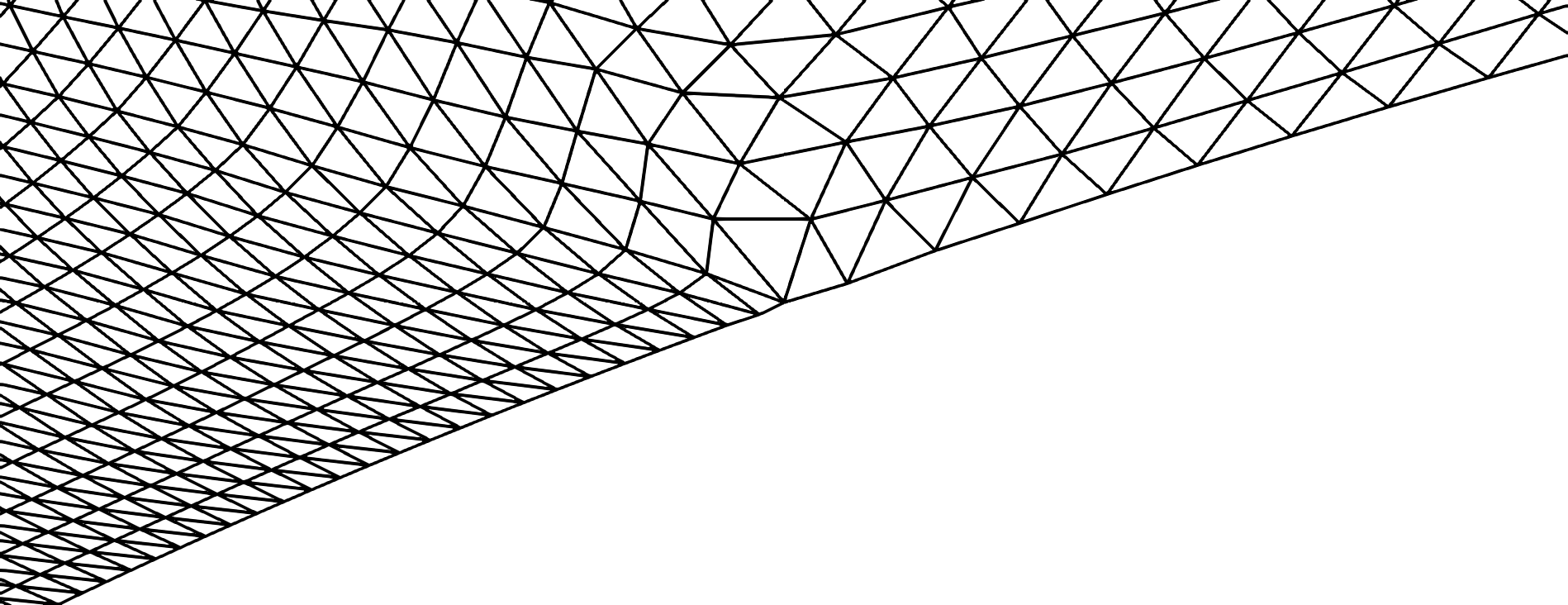}& \includegraphics[width=0.22\textwidth]{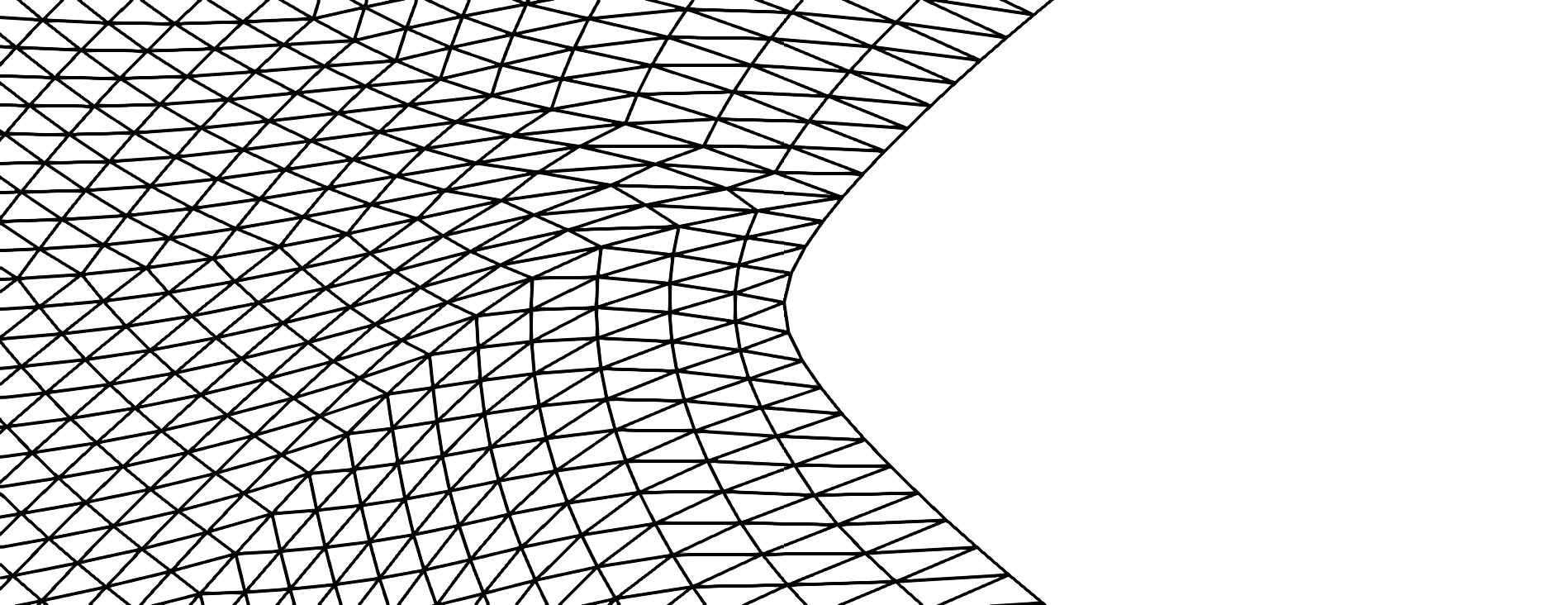}\\			
		\end{tabular}
	\end{center}
	\caption{Selected steps \{0,5,10,15,20,40\} of the optimization process show the large deformations occurring during the initial steps. The upper left corner and front tip are featured together with the profile of $\Gobs$. }
	\label{fig:2d-sequence}
\end{figure}

We provide in Figure \ref{fig:2d-sequence} a few of the shape iterates.
These demonstrate the large deformations which the optimization scheme provides in the first steps.
In the first column, shape 5 already shows the removal of the preexisting corners of $\Gobs$, as well as how the obstacle stretches across the direction of the flow.
This can be observed in detail in the next columns, where the elements across the upper left corner and the front tip are shown.
It is observed that large deformations are possible without provoking a failure of the numerical solver due to highly degenerate elements.
Moreover we wish to emphasize, besides solver tolerances, there are no quantities to tune, especially quantities which are designed to preserve mesh quality.

\subsection{Simulations in 3d}
For the 3d case a grid with \num{622592} tetrahedral elements was used. 
The surface of the obstacle consists of \num{12288} triangular elements.
The derivative $J'$ is scaled by a factor of $t=0.1$ for the computation of descent direction.
Moreover, for computational reasons, we utilize the Frobenius norm for \eqref{eq:geometricDirectionOfDescent} in this setting.
The outer loop in \Cref{alg:descent_alg} is again set to stop after 50 iterations.

 \begin{figure}[!htbp]
 	\centering
 	\includegraphics[width=0.8\textwidth]{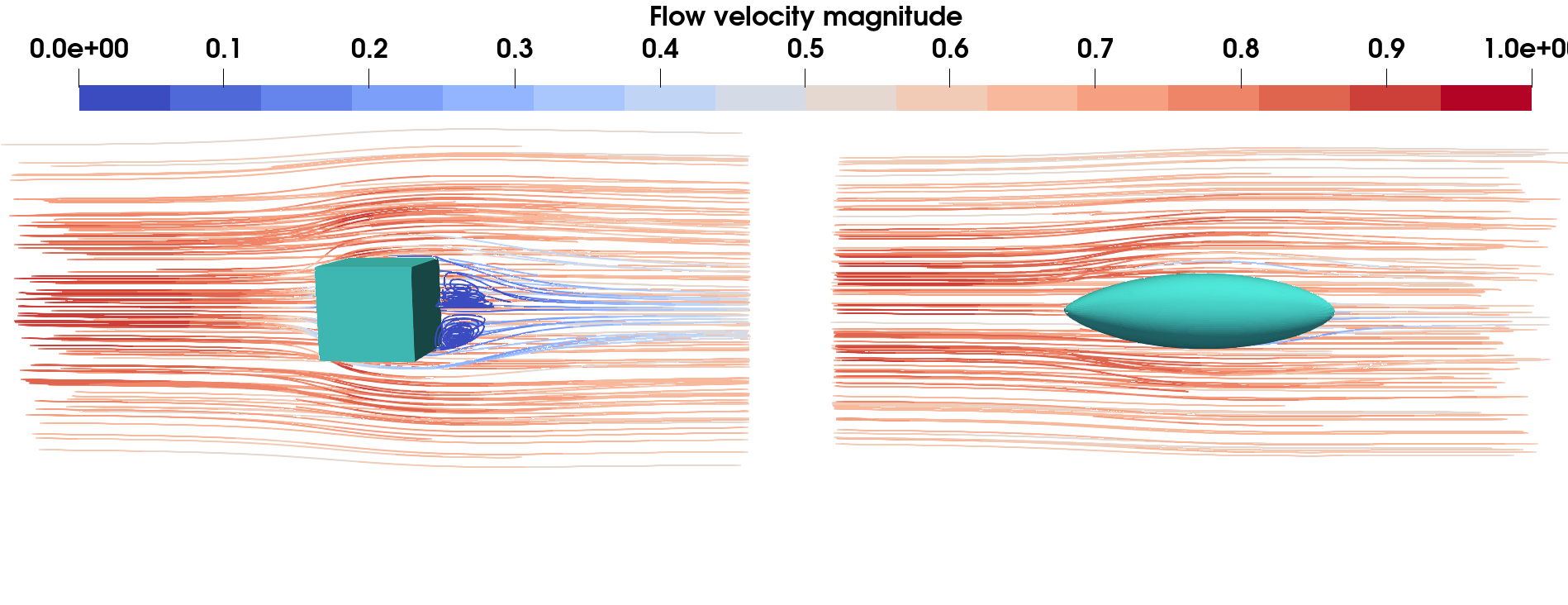}\\
 	\includegraphics[width=0.45\textwidth]{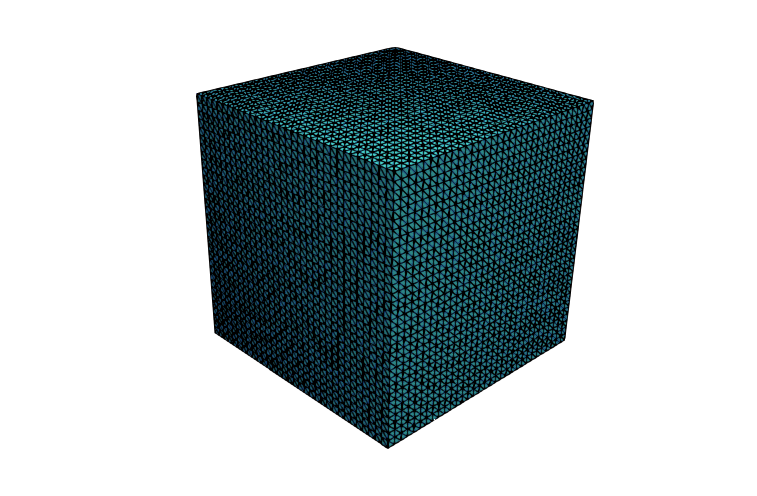}
 	\includegraphics[width=0.45\textwidth]{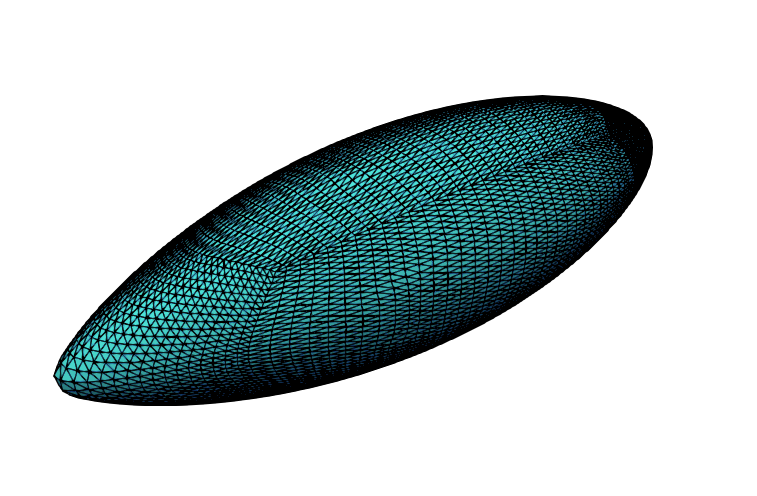}
 	\caption{3d results for a grid with \num{622592} tetrahedrons. The optimized obstacle is shown next to the reference. The streamlines around the object are shown (top), together with the grids (bottom).}
 	\label{fig:3dGrids&Streams}
 \end{figure}

The initial and final step of a 3d simulation are shown in \Cref{fig:3dGrids&Streams}. 
The streamlines show the regions on the initial grid where there is a disruption in the flow, leading to a higher energy dissipation.
As in the 2d case, the edges and corners of the box must be removed. 
This is shown in the lower row of \Cref{fig:3dGrids&Streams}, where the optimized geometry can be seen.
The triangular elements are shown over the surface of the object, which allows one to see the large deformations necessary for the removal of the geometric singularities.
Additionally, the creation of tips leads to a more optimal flow across $\Gobs$, as can be seen on the streamlines of the optimized geometry.

\begin{figure}[!htbp]
	\centering
	\includegraphics[width=0.7\textwidth]{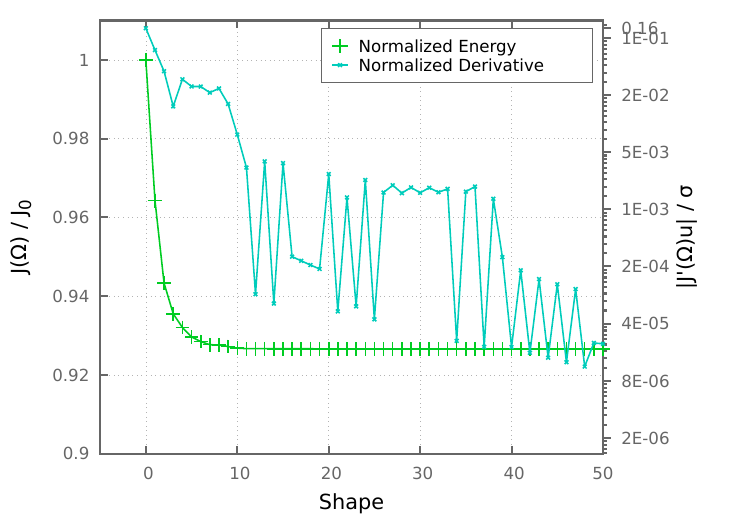}
	\caption{The objective function $J(\Omega)$, for a 3d simulation, divided by the initial value  - $J_0$ - over 50 shapes is shown, together with the convergence criterion \eqref{eq:ConvCriterion}.}
	\label{fig:PlotObjFunc3D}
\end{figure}

The results in \Cref{fig:3dGrids&Streams} can be linked to the plots given in \Cref{fig:PlotObjFunc3D}. 
As the preexisting singularities, i.e.~the edges and vertices, are removed, the dissipated energy decreases  pronouncedly. 
By step 10, an 8\% reduction has been achieved.
This is similar to the sequence shown in \Cref{fig:2d-sequence}, where in the early steps the obstacle is elongated parallel to the flow direction.
Together with the streamlines for the optimized geometry in \Cref{fig:3dGrids&Streams}, it can be seen how the more uniform flow around the surface results in a lower energy dissipation. 

\section{Comparison of descent directions in $\InfBanach$ and $\PBanach$}\label{sec:comp}
In this section a comparison between the approach given in \cite{mueller2022scalable} for vector fields in $\PBanach$ and $\InfBanach$ is performed.
The so-called $p$-Laplace relaxation scheme was proposed in \cite{deckelnick2021} in a shape optimization context and first applied in \cite{mueller2021}.
It is based on solving a relaxed problem in $\PBanach$ which is meant to approximate a steepest descent in $\InfBanach$.
The work \cite{Herbert23} discusses some of these approximations.
This relaxed formulation was inspired by \cite{ishii2005limits} which considers scalar functions.
Even moderate values of $p$ are found to be useful, however high values of $p$ are necessary to yield a good approximation to $\InfBanach$.

As in \cite{mueller2022scalable}, let the Lagrangian
\begin{equation}
	\begin{aligned}
		L^p(u,\lambda) :=\ &
		J'(\Omega)\, u + \frac{1}{p} \int_{\Omega} (Du : Du)^{p/2} \, dx + \sum_{i = 1}^{d+1} \mult{g,i} g_i(\Omega(\defor)),
	\end{aligned}
	\label{eq:PLagrangian}
\end{equation}
be used to obtain a highly nonlinear optimality system, which is solved for $\defor$. 
For moderate values of $p$, the solution is found by using increasing values of $p$, which compute the initial guess for the next increment, up to a given $p_{\mathrm{max}}$.
The latter is a caveat on itself, given the high computational requirements of successively solving a problem, which depends on an arbitrary increment to $p$. 
Additionally, as described in \cite[Sec.2.1]{mueller2022scalable}, care must be taken to prevent numerical problems associated to the second derivative of the Lagrangian, \eqref{eq:PLagrangian}.
The algorithm proposed in \cite{mueller2022scalable} is tailored in such a way that the largest possible deformations are allowed per step, without loss of solver convergence, which is a limiting factor for this $p$-Laplace algorithmic approach.
An in-depth discussion and high-performance computing results are presented in \cite{mueller2022scalable}. 

In order to compare both optimization schemes, the same energy dissipation problem \eqref{eq:EnergyDispObjectiveFunctional} is used.
The viscosity continues to be given by $\visc=0.02$.
A 2d computational mesh with four refinement levels and \num{70656} triangular elements is used for these simulations.
As in \cite{mueller2022scalable}, for the $p$-Laplace algorithm the maximum value of $p$ is set to $p_{\mathrm{max}}=4.8$. 
The ADMM-based optimization is configured as described in \Cref{sec:results}.
Similarly, the initial configuration features a box-shaped obstacle, thus it is necessary to remove and create geometric singularities. 

\begin{figure}[!htbp]
	\centering
	\includegraphics[width=0.6\textwidth]{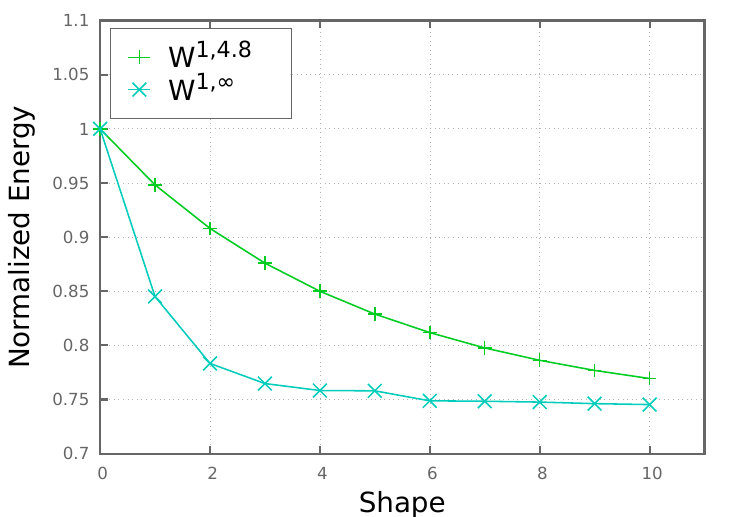}
	\caption{The dissipated energy for each approach, divided by their corresponding initial value - $J_0$- for the first 10 shapes is given. }
	\label{fig:ObjFunComp}
\end{figure}

Given that a highly viscous model is used, a large reduction of the dissipated energy is not necessarily related to the appearance of sharp geometric singularities over the surface of obstacle, $\Gobs$, but to the deformation of the obstacle into the well-known prolate spheroid \cite{pironneau1973optimum,pironneau1974optimum}. 
As shown in \Cref{fig:2d-sequence}, for the ADMM this occurs within the initial 10 steps of \Cref{alg:descent_alg}.

The reduction of the objective function for the first 10 steps is illustrated in \Cref{fig:ObjFunComp}.
The values of $J$ are divided by the initial value $J_0$, with the purpose of allowing for a correct visualization of the differences between the two methods. 
The ADMM-based algorithm with $\defor\in\InfBanach$ allows for a 15\% objective function reduction within the first shape optimization step.
On the contrary, the $p$-Laplacian approach reaches this reduction level only after about four shape iterates. 
Around step 6, the $\InfBanach$ approach allows for a 25\% decrease on the initial dissipated energy. 
The relaxed approach required more than 10 steps to achieve this level of reduction, as seen in the figure.

\begin{figure}[!htbp]
	\centering
	\includegraphics[width=0.55\textwidth]{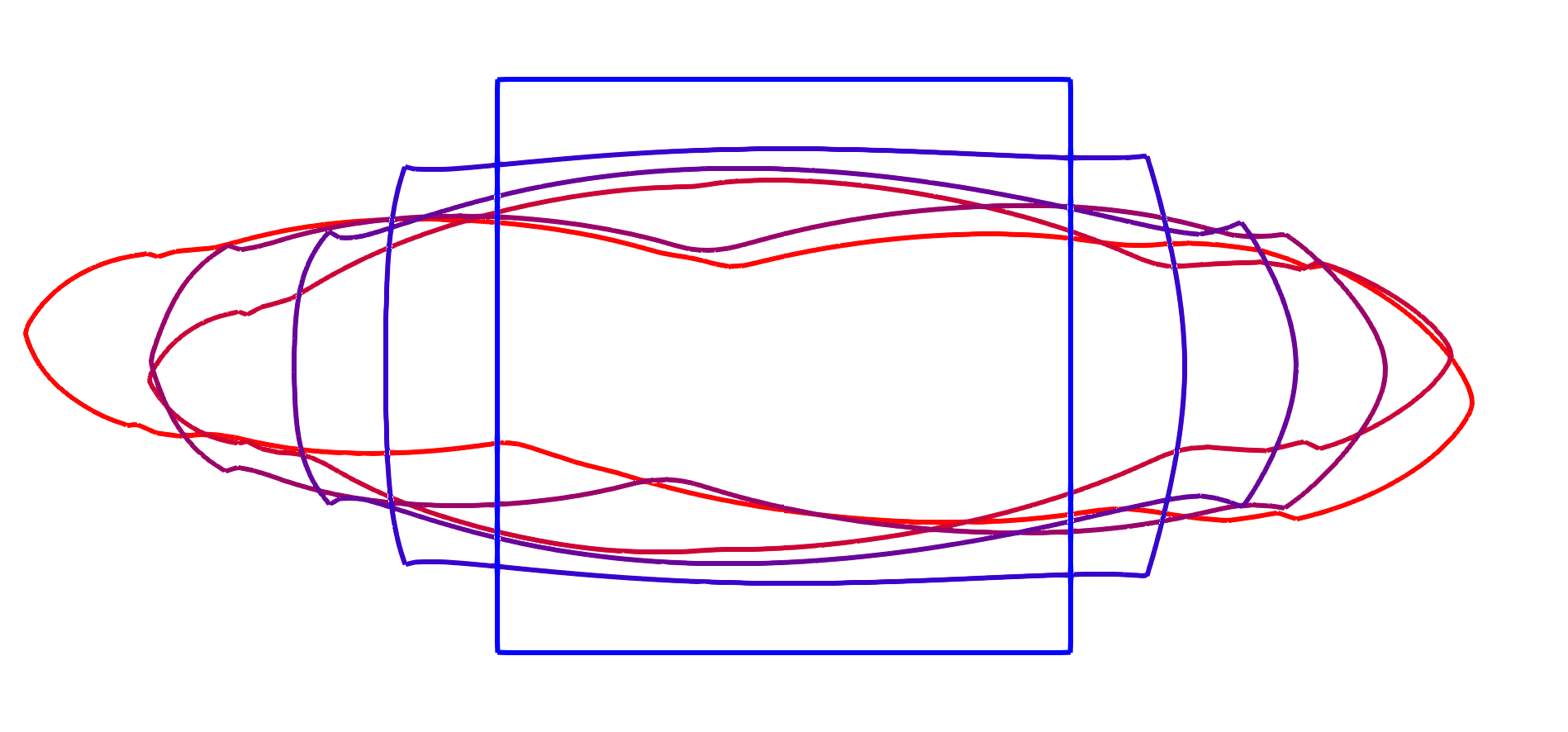}\hspace{-1.7cm}
	\includegraphics[width=0.55\textwidth]{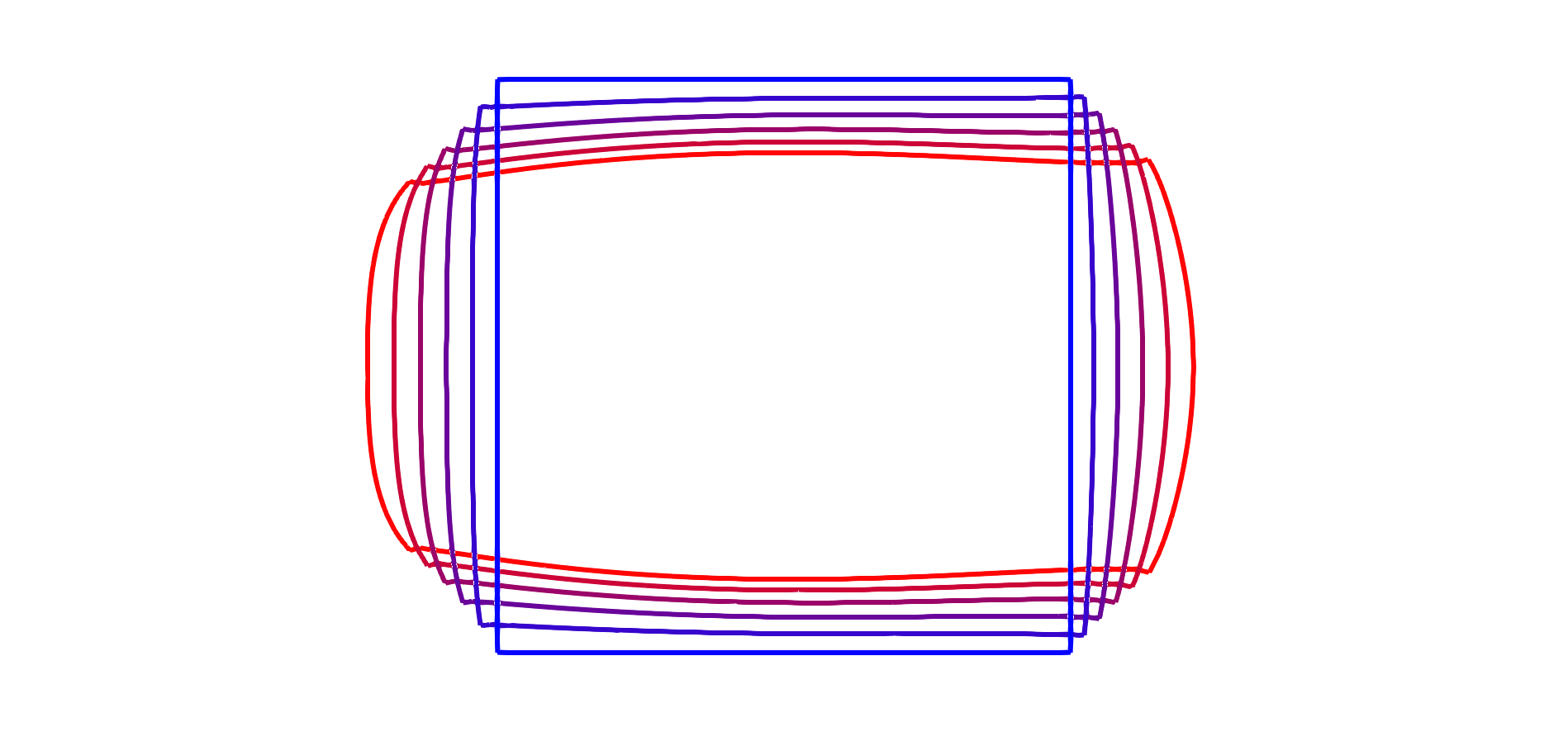}
	\caption{The initial and first five shape iterates for the $\InfBanach$ (left) and the $\PBanach$ (right) shape optimization schemes. The initial configuration $\Omega_0$ (blue) is deformed by applying the computed deformation field until $\Omega_5$ (red) is obtained.}
	\label{fig:BoxShapesComp}
\end{figure}

When observing the shapes rather than the energy, one might attribute the large differences of the energies between the methods to the aggressive deformations which appear from the $\InfBanach$ method.
A comparison of the first five shapes for the $\PBanach$ and $\InfBanach$ are presented in \Cref{fig:BoxShapesComp}.
It may be seen that, while we expect that they will eventually provide the same shape, their paths to becoming an optimized shape are rather different.
The deformation fields with the $\PBanach$ method gently deform $\Gobs$, slowly stretching it out.
As the geometry is updated, the edges perpendicular to the flow are elongated and acquire a rounded profile.
Using $\InfBanach$ for the descent it immediately generate this oblong and elongated shape which is seen as the fifth shape for the $p$-Laplace approach.
The appearance of the tips can already be observed in step 2, together with a large elongation parallel to the flow direction.
By steps 4 and 5, the tips have been created and the geometric singularities, i.e.~the corners of the box, have mostly been removed. 

\begin{figure}[!htbp]
	\centering
	\hspace{1cm}\includegraphics[width=0.45\textwidth]{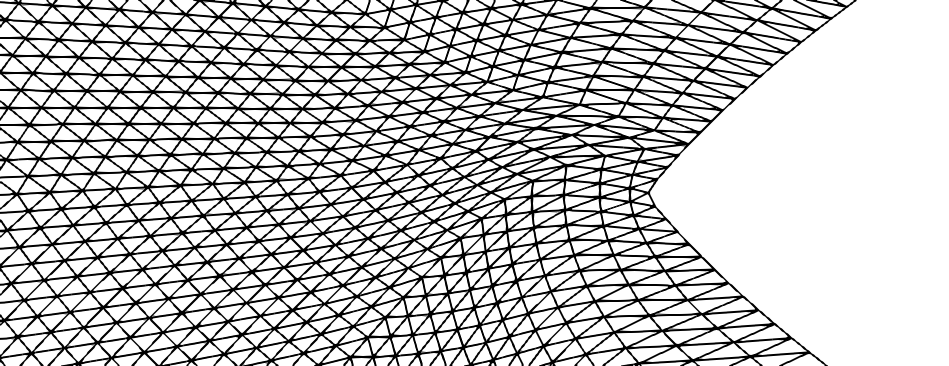}
	\includegraphics[width=0.45\textwidth]{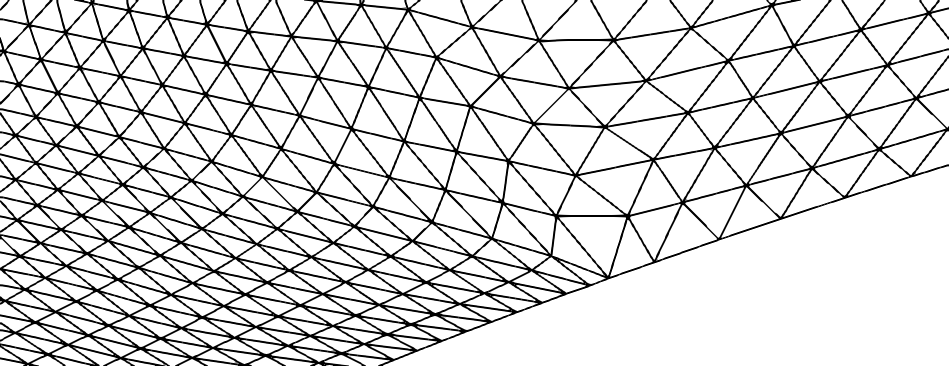}\\
	\vspace{0.3cm}
	\hspace{1cm}\includegraphics[width=0.45\textwidth]{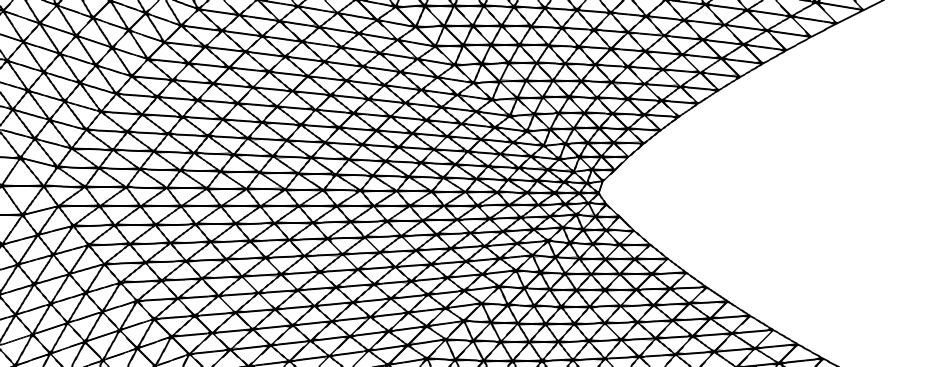}
	\includegraphics[width=0.45\textwidth]{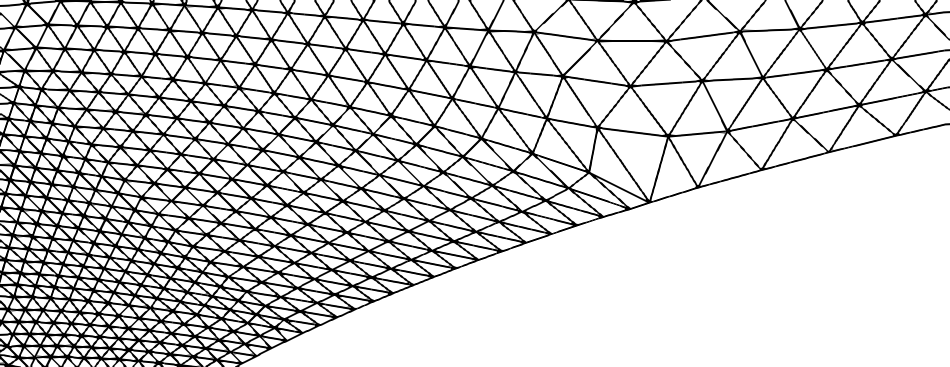}
	\caption{Computational meshes, corresponding to the optimized shapes after 50 steps, for a $\InfBanach$ (top) and $\PBanach$ (bottom) optimization schemes. The front tips (left) and upper left corners (right) are shown.}
	\label{fig:FinalGridsComp}
\end{figure}
However, these large deformations during the early stages of the optimization process may have a negative impact on mesh quality.
In the context of the finite element discretization, it is of interest that the optimized shape retains mesh quality, so that remeshing can be avoided.
For instance, the computational grid must allow for the efficient and accurate solution of the state equation and the computation of the objective function.
A particular concern is the overlapping of the grid, which is prevented by the $\InfBanach$ method, and degeneracy of triangles (in 2d) which is not guaranteed.

The grids of the optimized shape are compared in \Cref{fig:FinalGridsComp}. 
Two regions around $\Gobs$ are studied, the tip frontal to the flow direction and the section of the grid where the upper left corner of the box has been smoothed out. 
Given that the $\InfBanach$ vector fields largely deform these regions in the initial steps, the effects can be seen in the triangles over the surface of the obstacle. 
Compared to the $p$-Laplace descent directions, the elements have moved around the $\Gobs$.
The surrounding elements are more evenly distributed in the relaxed case, compared to the uneven pattern at the front tip of the $\InfBanach$.
Nevertheless, the large deformations enabled by the ADMM do not result in degenerate elements, particularly on the critical areas shown.

\begin{table}[!htbp]
	\caption{Edge quality data for different refinement levels.
			The edge length ratio is used as a quality metric, i.e. the ratio of the longest to the shortest edge over $\Gobs$.}
	\begin{center}
		\begin{tabular}{||c:cc:cc:cc||}
			\multicolumn{1}{c}{}&\multicolumn{2}{c}{3-refs}&\multicolumn{2}{c}{4-refs}&\multicolumn{2}{c}{5-refs}\\
			\hline\hline
			\multicolumn{1}{c}{Triangles in $\Omega$ }&\multicolumn{2}{c}{\num{17664}}&\multicolumn{2}{c}{\num{70656}}&\multicolumn{2}{c}{\num{282624}}\\
			\multicolumn{1}{c}{Edges in $\Gobs$ }&\multicolumn{2}{c}{\num{128}}&\multicolumn{2}{c}{\num{256}}&\multicolumn{2}{c}{\num{512}}\\
			\hline
			Step&ADMM&PLAP&ADMM&PLAP&ADMM&PLAP\\
			\hline
			0&1.00 &1.00& 1.00 & 1.00& 1.00 & 1.00\\
			1&1.89 &1.20& 1.91 & 1.23& 1.96 & 1.28\\
			2&3.62 &1.42& 3.67 & 1.49& 3.61 & 1.59\\
			3&5.58 &1.64& 5.27 & 1.76& 6.45 & 1.91\\
			4&10.54&1.87& 10.57& 2.03& 11.01& 2.24\\
			5&10.52&2.09& 10.53& 2.29& 9.71 & 2.56\\
			\vdots&&&&&&\\
			47& 13.49 & 8.33 & 14.96& 9.56 & 13.69& 8.86\\
			48& 13.50 & 8.51 & 15.05& 9.89 & 13.64& 8.98\\
			49& 13.50 & 8.54 & 15.03& 9.73 & 13.64& 9.09\\
			50& 13.50 & 8.56 & 15.12& 10.17& 13.60& 9.30\\
			\hline\hline
		\end{tabular}
	\end{center}
	\label{tab:SurfaceEdgeDataBox}
\end{table}
The differences in the resulting grids, across the whole optimization process, can be quantified by the effect of the deformations on the elements that conform the obstacle surface $\Gobs$.
In the 2d computational grid, the surface of the obstacle is formed by the edges of the adjacent triangles.
This is the region of the domain where the geometric singularities are to be created and removed, thus it is one of the regions that undergoes large deformations.
In this sense, it is of interest to study the effect that these deformations have on the quality of the surface $\Gobs$.

In \Cref{tab:SurfaceEdgeDataBox}, tests were carried out for several levels of refinements, up to \num{282624} elements in $\Omega$ and 512 edges in $\Gobs$. 
The simulations were configured to run for a total of 50 steps.
The edge length ratio, that is ratio between the longest and the shortest edge of $\Gobs$, is shown for the first 5 and last 4 shape iterates.
A comparison is made between the grids resulting from simulations based on $\InfBanach$ and $\PBanach$ descent directions.

The evenly distributed grid allows for the edge length ratio to be equal to a unit for the initial grid $\Omega_0$. 
In line with the previously mentioned aggressive objective function reduction and with the observed large deformations in the initial steps, the mesh quality data in \Cref{tab:SurfaceEdgeDataBox} shows how the ratio increases faster for the $\InfBanach$ than for the $p$-Laplace method. 
This trend is repeated across all refinement levels.
At the end of the simulation the quality metric for the edges is higher for the $\InfBanach$ case.
Nonetheless, the ADMM algorithm allows for larger deformations without provoking negative effects on the iterative solvers used for the system of equations.
While such large deformations would most likely not be reachable by the $p$-Laplace approach or would have a negative impact on the convergence of the solvers used to compute $\defor$, using the ADMM it is possible to promote the preservation of mesh quality while allowing for a very aggressive shape optimization.

\section{High-performance computing results}\label{sec:hpc}
In this section we present a parallel scalability study for weak scaling in 2d simulations. 
The results were computed using the supercomputer Lise at HLRN. 
This distributed-memory system consists of a total of 1270 compute nodes, each with a dual-socket architecture of 48 cores per CPU at a maximum frequency of 2.30GHz and 384GBs of RAM. 
The initial core count is set to 48 cores, all assigned within one socket. 
This study was carried out for up to \num{3072} cores, distributed across 32 nodes at full usage.

\begin{figure}[!htbp]
	\centering
	\begin{subfigure}[t]{\textwidth}
		\centering
		\scalebox{0.46}{\includegraphics[width=0.999\textwidth]{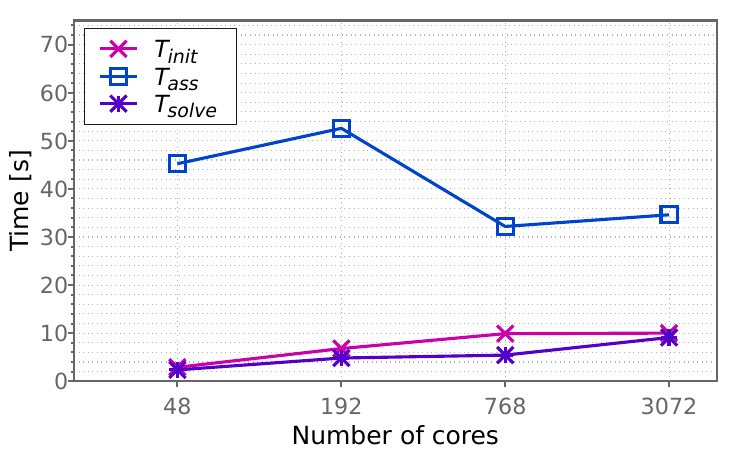}}
		\scalebox{0.46}{\includegraphics[width=0.999\textwidth]{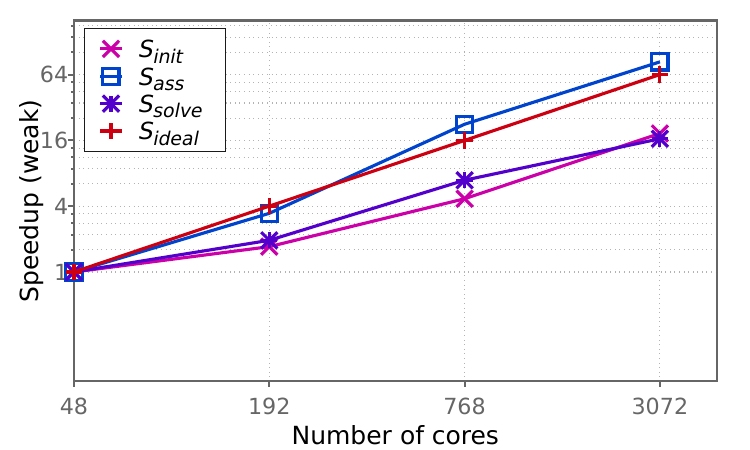}}	
		\caption{Weak scalability: time measurements and speedup relative to 48 cores.}
		\label{fig:TimeAndSpeedup}
	\end{subfigure}	
	\vskip\baselineskip
	\centering
	\begin{subfigure}[t]{\textwidth}
		\centering
		\resizebox{\columnwidth}{!}{
		\begin{tabular}{rrrccc}
			\centering
			\begin{tabular}{llllcccc}
				\toprule
				Procs  & Refs  & Number of& DoFs& AVG ADMM & AVG Newton &AVG Lin.Its.& AVG Time\\
				& & Elements& & Its. per Shape&Its. per Shape&  per Shape &per Shape[s]\\
				\midrule
				\num{48}  & 3 & \num{17664}   & \num{18016}   & 31& 166& \num{5267} & 59\\
				\num{192} & 4 & \num{70656}   & \num{71360}   & 37& 197& \num{6825} & 76\\
				\num{768} & 5 & \num{282624}  & \num{284032}  & 25& 139& \num{5068} & 56\\
				\num{3072}& 6 & \num{1130496} & \num{1133312} & 25& 141& \num{5463} & 65\\
				\bottomrule
			\end{tabular}
		\end{tabular}
		}
		\caption{Average iteration counts for the weak scaling results}
		\label{tab:iter_counts}
	\end{subfigure}
	\caption{
	Parallel scalability results for the first five optimization steps in a 2d simulation.
	Accumulated time and gained speedup relative to 48 cores are shown for the linearization solved within the minimization with respect to $\defor$. 
	Average iteration counts are given for the ADMM routine and Newton's method used therein.}
	\label{fig:scaling}
\end{figure}

A 2d mesh with \num{17664} triangular elements was used as starting point for these measurements. 
The number of elements is increased by a factor of four with each level of refinement.
Accumulated wallclock times and speedup relative to 48 cores are provided for up to a million elements for the first five optimization steps.
Measurements of the time for the assembly, $(T_{\mathrm{ass}})$,  of the linearized system; the initialization of the preconditioner, $(T_{\mathrm{init}})$; and the application of the linear solver, $(T_{\mathrm{solve}})$, for each step within Newton's method are presented. 
Specifically, line 5 of \Cref{alg:ADMM} is measured, because it is the most computationally relevant aspect of the routine.
It is where the majority of the computational effort is spent.
Since it is not the focus of this article, the time required for the solution of the state and adjoint equation systems are excluded from the displayed measurements. 

In addition to the mentioned wallclock times, average iteration counts of the ADMM routine are analyzed, together with an estimate of the average time per shape of the whole outer loop where shape iterates are created. 
For the ADMM, an iteration represents one full cycle of the loop in lines 3-18 of \Cref{alg:ADMM}.
From this loop, the average calls to Newton's method are given, including the linear iterations required to find $\defor$ in line 5 of \Cref{alg:ADMM}.
The time-per-shape measurement is comprised of all the calls to the routine \Cref{alg:ADMM} within an iteration of the outer loop in \Cref{alg:descent_alg}, excepting for the previously mentioned fluid solver quantities.

Further to the algorithmic details given in \Cref{sec:optim}, we now describe the implemented linear algebra solvers.  
The linearization within \eqref{eq:optimality_system_defor} is solved using a BiCGStab method preconditioned by a geometric multigrid.
The preconditioner uses 3 pre- and post-smoothing Gauss-Seidel steps with a V-cycle configuration.
The linear solver is configured to reach a relative and absolute reduction of $10^{-10}$ and $10^{-12}$, respectively. 
Moreover, the coarsest grid is solved in a single core via an LU factorization. 
For the solution of \eqref{eq:optimality_sub-system_q}, a BiCGStab method is used with a Jacobi preconditioner. 
Given the structure of this problem, the system can be solved in a single iteration.
 
The weak scalability study is presented in \Cref{fig:scaling}, where scalable results are seen for up to \num{3072} cores. 
The accumulated times show a slight performance drop for a 192 core distribution in the assembly phase. 
Also note that there is a coinciding increase in the average calls to the ADMM routine.  
Correspondingly, there is an increase in the average time required to generate the shape iterates.
For the latter core counts, the average computational cost remains roughly equal, as observed in the iteration counts table.
The number of degrees of freedom (DoFs) refers to the size of the system used to compute $\defor$, which is equal to $d$ times the number of vertices of the mesh. 
This equation system increases in size by two orders of magnitude at six levels of refinement. 
Nevertheless, the measurements in \Cref{fig:scaling} show that the computational load is bound throughout all refinement levels.

\section{Conclusion}\label{sec:concl}
In this paper, we presented a shape optimization methodology that uses the ADMM to obtain descent directions in $\InfBanach$.
The advantages of this approach were illustrated through 2d and 3d results, which show how this methodology allows for large deformations in the early steps of the optimization process.
Our formulation introduced nonlinear geometric constraints, necessary to avoid trivial solutions. 
The method to handle these constraints was taken from \cite{mueller2022scalable} and it allows for fulfillment of the geometry constraints on each generated shape.

Weak parallel scalability results for up to \num{3072} processes in a distributed-memory system were presented. 
The use of appropriate preconditioners, with mesh-independent convergence properties, were utilized to achieve an scalable implementation.
This is accentuated within the large number of linear iterations required for one shape optimization step. 
In our studies it was seen that the geometric multigrid method bounds the average computational workload, even under large increments of the number of DoFs. 
Additionally, a comparison of the early optimization steps was performed between the approach here proposed and the $p$-Laplace method as formulated in \cite{mueller2022scalable}.
This showed that larger deformations are encouraged with a $\InfBanach$ optimization scheme.
Per shape, it is seen that the objective function decreases faster for the $\InfBanach$ method. 
A mesh quality study shows that the $\InfBanach$ approach results in a not-meaningful decrease in quality of the line elements which compose the surface of the obstacle, as compared to the $\PBanach$ algorithm. 
Taking into account that large deformations are needed to take the initial geometry to the optimized domain, a decrease in mesh quality is expected. 
However, in industrial applications, by expert-knowledge, the initial mesh is likely to be closer to the sought optimizer.   

In this work a fixed step-size ADMM is used, unlike in \cite{deckelnick23} where the variable step-size version is implemented. 
It is worth mentioning that the methodology here presented is not bound to the ADMM, i.e. another method could be used. 
However, the area of shape optimization in $\InfBanach$, and the algorithms required to compute adequate descent directions, is an ongoing field of research. 

We suggest that our proposed technique could be considered for more complex geometries. 
The topic of mesh quality, in terms of its impact on the state variable, is to be explored in further work to widen our understanding of the advantages and disadvantages of different shape optimization methodologies.

 \subsection*{Acknowledgements}
 {	
 	\noindent
 	The current work is part of the research training group ``Simulation-Based Design Optimization of Dynamic Systems Under Uncertainties'' (SENSUS) funded by the state of Hamburg under the aegis of the Landesforschungs\-förderungs-Project LFF-GK11.
 	
 	\noindent
 	Computing time on the national supercomputer Lise at the North German Supercomputing Alliance (HLRN) under the project hhm00006 (An optimal shape matters) is gratefully acknowledged.
 }

\printbibliography
\end{document}